\documentclass[a4paper,10pt]{article}
\usepackage[utf8]{inputenc}
\usepackage{graphicx}       
\usepackage{longtable} 
\usepackage{color}
\usepackage{amssymb,amsmath}
\usepackage{float}
\usepackage{enumitem}
\usepackage{subfigure}
\usepackage{units}

\usepackage{authblk}

\newcommand{\bs}[1]{\boldsymbol{#1}}
\newcommand{\R}{\mathbb{R}}

\title{Hybrid models for simulating blood flow in microvascular networks}
\author{E. Vidotto$^1$, T. Koch$^2$, T. K\"oppl$^1$, R. Helmig$^2$, B. Wohlmuth$^1$}

\date{\vspace{-5ex}}

\begin{document}

\maketitle

\footnotetext[1]{Department of Mathematics, University of Technology Munich, 
Boltzmannstr. 3, 85748 Garching bei M\"unchen, Germany, vidotto@ma.tum.de, koepplto@ma.tum.de, wohlmuth@ma.tum.de}

\footnotetext[2]{Department of Hydromechanics and Modelling of Hydrosystems, University of Stuttgart, Pfaffenwaldring 61, 70569 Stuttgart, Germany
timo.koch@iws.uni-stuttgart.de, rainer.helmig@iws.uni-stuttgart.de}

\begin{abstract}
In this paper, we are concerned with the simulation of blood flow in microvascular networks and the
surrounding tissue. To reduce the computational complexity of this issue, the network structures are modeled by a one-dimensional
graph, whose location in space is determined by the centerlines of the three-dimensional vessels. The surrounding tissue is considered as a
homogeneous porous medium. Darcy's equation is used to simulate flow in the extra-vascular space, where the mass exchange with the blood
vessels is accounted for by means of line source terms. However, this model reduction approach still causes high computational costs,
in particular, when larger parts of an organ have to be simulated. This observation motivates the consideration of a further model reduction step. Thereby,
we homogenize the fine scale structures of the microvascular networks resulting in a new hybrid approach modeling the fine scale structures as a
heterogeneous porous medium and the flow in the larger vessels by one-dimensional flow equations. 
Both modeling approaches are compared with respect to mass fluxes
and averaged pressures. The simulations have been performed on a microvascular network that has been extracted from the cortex of a rat brain.
\end{abstract}

\section{Introduction}
Modeling of blood flow and transport at the level of microcirculation is an interesting subfield in
biomedical engineering. A reliable computational model for the microcirculation of the human body
would enable physicians and pharmacists to obtain better insight
into the oxygen supply of cells, the waste removal from the interstitial space \cite{holter2017interstitial}
and further important biological processes without
the need to perform expensive and risky experiments \cite{boas2008vascular,DAngelo,pries2014making}. Besides, such models open up the
possibility to study the impact of diseases like Alzheimer's \cite{iturria2016early} and to improve medicines and therapies for cancer treatment
\cite{cattaneo2014computational,hsu1989green,mahjoob2009analytical,nabil2016computational}.
Well known therapies for cancer treatment are, e.g., hyperthermia \cite{petryk2013magnetic,waterman1991blood}
or the injection of therapeutic agents preventing the vascularization of cancer tissue \cite{stacker2018emerging}.
In order to increase the efficiency of the mentioned cancer therapies, it is crucial to focus the therapeutic agents or the heat
on the cancerous part of an affected organ while maintaining the rest of the tissue. In order to be able to
simulate the distribution of heat and the concentration of therapeutic agents with an adequate accuracy, the hierarchical structure of the
vascular system supplying the considered organ has to be taken into account.

However, even for vascular systems contained in small
volumes covering just a few cubic millimetres, it is a challenging task
to simulate flow in the microvascular networks and the interstitial space, since the flow in a complex network structure consisting of
thousands of vessels \cite{formaggia2010cardiovascular} (Chap. 1) is coupled with the flow in the surrounding 
tissue matrix \cite{dewhirst2017transport,pries1990blood}.

As a consequence, several model reduction techniques have been developed for simulating flow through an entire
organ or parts of an organ. A widespread technique in this field are homogenization techniques. Thereby brain tissue and the vascular
system are considered as two different continua modeled as porous media
\cite{el2015multi,penta2015multiscale,roose2012multiscale,shipley2010multiscale} and the
pressure and velocity field are computed by means of Darcy's equation \cite{helmig1997multiphase,khaled2003role}.
Using this approach, the computational
effort and the data volume are significantly reduced. On the other hand, both the pressure and the velocity field can only be described in
an averaged sense, neglecting the exact structure of the vascular system. As a result, the distribution of therapeutic agents or heat might not be computed
with sufficient accuracy. A further challenge is to determine suitable permeability tensors for the Darcy equation reflecting the
structure of the vascular trees in an averaged way. Thereby brain tissue and the vascular system are considered as 
two different continua modeled as porous media
\cite{el2015multi,penta2015multiscale,roose2012multiscale,shipley2010multiscale} and the
pressure and velocity field are computed by means of Darcy's equation \cite{helmig1997multiphase,khaled2003role,stoverud2012modeling,penta2015multiscale}.\\
Another way to decrease the computational complexity, is to describe the vascular networks
by means of one-dimensional (1D) flow models \cite{d2008coupling,tang2018numerical,nabil2016computational,formaggia1999multiscale,
quarteroni2004mathematical,gjerde2018splitting}, while the surrounding tissue is considered as a three-dimensional porous medium. By this, an expensive
meshing of a three-dimensional (3D) blood vessel system is avoided and at the same time the hierarchical structure of the vascular system is maintained.
However, elaborate concepts for coupling one-dimensional flow equations with a
flow equation (Darcy equation) for a three-dimensional flow problem \cite{d2012finite,koppl2016local,koppl2018mathematical} have to be developed.
This is done by constructing suitable source terms for both the
Darcy equation in 3D and the 1D flow problems. On closer examination, one notes that the source term in 3D consists of Dirac measures concentrated
on the centerline of the vessels or the corresponding vessel walls, inducing a certain roughness into the solution of the 3D Darcy equation.
Due to possible singularities or kinks occurring
in the solution sufficiently fine meshes are required to obtain an accurate numerical solution.

Here, we propose a new hybrid approach which preserves the advantages of reduced order models, while preserving a sufficiently high accuracy. The idea is to
model the larger vessels by 1D flow models, whereas the capillaries and tissue are considered as porous continua, as it is done in the homogenization approach.
We apply this hybrid modeling approach to a microvascular system filling a volume of $1.134\;\unit{mm} \times 1.134\;\unit{mm} \times 2.268\;\unit{mm}$
(see Fig. \ref{fig:extracted_network}). This microvascular system was taken from the brain of a rat and the corresponding data was generated by
the group of B. Weber (University of Zurich) \cite{reichold2009vascular}. Within this network, one can  find several 
penetrating arterioles and venules connecting this vascular subsystem to the macrocirculation of the rat's brain, see Fig. \ref{fig:extracted_network}.
Our motivation to choose this approach is that by resolving the largest vessels in this volume, the main hierarchy of the vessel system is 
incorporated into the model.
Moreover, using a homogenized double continuum model for the capillaries and the tissue, we have a less complex
model which will certainly yield some speed-up, if
several units of cubes of a few cubic millimetres are aggregated to simulate larger parts of an organ. 
Compared to other works \cite{el2018investigating,peyrounette2018multiscale,shipley2019hybrid}
investigating hybrid approaches for microvascular systems, 
we consider in this work the capillaries and the tissue as two coupled porous media, i.e., a double continuum approach \cite{erbertseder2012coupled}. 
In addition to that we discuss alternative coupling conditions between the 1D vessels and the 3D continuum for the capillaries. To
validate our \emph{hybrid 3D-1D} modeling approach, a comparison with a \emph{fully-discrete 3D-1D} 
model \cite{cattaneo2014computational,d2012finite,d2008coupling,notaro2016mixed,koppl2018mathematical}
is performed.

The rest of this work is structured as follows. In Section~\ref{sec:ProblemSetting}, we describe the data set that is used for our
numerical tests. In addition to
that the basic modeling assumptions are discussed. The model problems and the corresponding numerical 
discretization methods are outlined in Section \ref{sec:Models}.
Section~\ref{sec:numerical_test} contains some numerical tests and a discussion of the results obtained from both models. 
Finally in Section \ref{sec:Conclusion}, we summarize the key features of our work and give a short outlook.

\section{Problem setting and main simplifications}
\label{sec:ProblemSetting}
\begin{figure}[h!]
\begin{center}
\includegraphics[width=0.6\textwidth]{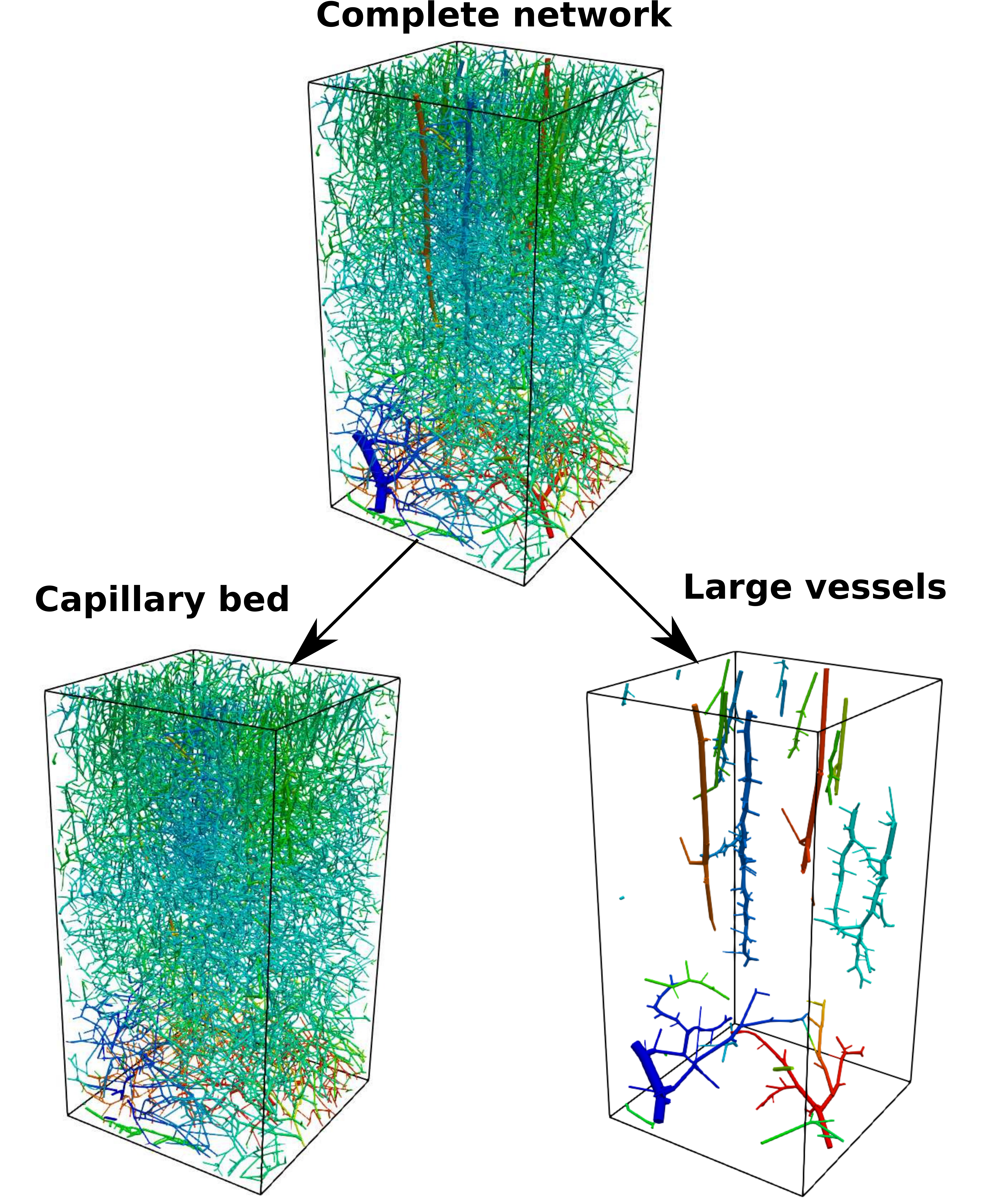}
\end{center}
\caption{\label{fig:extracted_network} A microvascular network extracted from a rat brain (top), filling 
a volume of $1\;\unit{mm} \times 1\;\unit{mm} \times 2\;\unit{mm}$.
A subdivision based on the vessel radius can be seen at the bottom of the figure.
The network at the bottom left, is the capillary bed of the microvascular network (all vessels whose radii are smaller than a given threshold, 
here $7\;\unit{\mu m}$).
The motivation for such threshold is postponed to Section~\ref{sec:R_T}. At the bottom right, arterioles and venules (whose radii are larger 
than the given threshold) constitute the network of larger vessels.}
\end{figure}

In order to illustrate the performance of our modeling approaches, we consider the microvascular network shown in Fig. \ref{fig:extracted_network}.
To obtain the data for this microvascular network, firstly an anaesthetized rat was perfused with phosphate buffered saline and barium sulphate. 
In a second step, the brain was removed and a sample from the cortex was taken. 

Using this sample, an X-ray tomographic microscopy of the sample was performed. 
Based on the resulting image, vessel midlines
were extracted, and the 3D coordinates of branches, kinks as well as start and end points were determined. We denote these points as network nodes.
Furthermore, curved vessel midlines were replaced by straight edges and to each edge, a mean vessel radius was assigned. Besides the different radii
also the connectivity of network nodes was determined, i.e.,
it was reported, which network nodes are connected with each other. At the network nodes that are adjacent to the boundaries of the extracted sample,
physiologically meaningful pressure boundary conditions are provided \cite{reichold2009vascular}.

Analyzing the given data, it turns out that the radii range from $1.6\;\unit{\mu m}$ to $28.2\;\unit{\mu m}$. Within the network some 
larger penetrating arterioles and venules
can be identified at the top of the sample. At the bottom of the sample a larger venule and arteriole are leaving or entering the cuboid domain. 
The rest of the microvascular
network consists of tiny capillaries. Zooming into the original data set, we detected some dead ends of the arterioles and venules 
(network nodes that are not connected to the
capillaries). Since we want to perform a precise comparison between the fully discrete and the hybrid modeling approach, the edges 
associated with dead ends are removed from
the data set such that in the resulting microvascular network all the network nodes are connected.
Further modeling assumptions that are used in the remaining sections, read as follows:
\ \\
\begin{itemize}
 \item[(A1)] The non-Newtonian flow behavior of blood is modeled in a simplified way using an algebraic relationship.
 \ \\
 \item[(A2)] The density of blood is constant.
 \ \\
 \item[(A3)] The influence of gravity is neglected.
 \ \\
 \item[(A4)] The inertial effects can be neglected.
 \ \\
 \item[(A5)] The pulsatility can be neglected.
 \ \\
 \item[(A6)] The walls of the larger vessels are considered impermeable, i.e., no flow occurs across their walls.
\end{itemize}
\ \\
These assumptions are motivated by the following considerations: It is well known that blood consists beside plasma also of blood cells. 
Blood plasma itself is a mixture of water and
ions, while the blood cells can be separated into the following groups: red blood cells, white blood cells 
and blood platelets \cite{formaggia2010cardiovascular} (Chap. 1).
In particular the red blood cells determine the flow behavior of blood significantly \cite{fahraeus1931viscosity}. 
Moving through a microvascular network the red blood cells
have to deform more and more as the vessel diameters become smaller than the diameter of the red blood cells. This results in a varying viscosity 
for the individual blood vessels
and therefore this feature is recorded in (A1)~\cite{gregorio2018computational}. Thereby, we adopt for an edge of diameter $D$ 
the following formula for the in vivo viscosity
$\mu_{\mathrm{bl}}\;\left[ \unit{Pa \cdot s} \right]$ \cite{pries1996biophysical}:
\begin{equation}
\label{eq:viscosity}
\mu_{\mathrm{bl}}\left( D \right) = \mu_{\mathrm{p}} \left( 1+ \left( \mu_{0.45} -1 \right)
\frac{\left(1-H \right)^C-1}{\left(1-0.45 \right)^C-1} \cdot \left( \frac{D}{D-1.1} \right)^2 \right)
\cdot \left( \frac{D}{D-1.1} \right)^2.
\end{equation}
Please note that in \eqref{eq:viscosity} the diameter $D$ is dimensionless. The physical diameter $d\;\left[\unit{\mu m}\right]$ has to be 
divided by $1.0\;\unit{\mu m}$ to obtain $D$.
The viscosity of the blood plasma is denoted
by $\mu_{\mathrm{p}}\;\left[ \unit{Pa \cdot s} \right]$, and $H$ stands
for the discharge hematocrit which is defined by the ratio between the volume of the red blood cells and the total blood volume. The
apparent viscosity $\mu_{0.45}$ for a discharge hematocrit of $0.45$ is calculated by:
$$
\mu_{0.45} = 6.0 \exp\left(-0.085 \cdot D \right) + 3.2 - 2.44 \exp\left(-0.06 \cdot D^{0.645} \right)
$$
and $C$ is a coefficient determining the influence of $H$ on $\mu_{\mathrm{bl}}$:
$$
C = \left( 0.8 + \exp\left(-0.075 \cdot D \right) \right) \left(-1 + \frac{1}{1 + 10^{-11} D^{12}} \right)
+ \frac{1}{1 + 10^{-11} D^{12}}.
$$
Please note that this function for the in vivo viscosity holds actually for human blood.
However we are not aware of such a function for rat blood, and therefore we assume that this function can be employed for rat blood as well.
For simplicity, we assume in (A2) that blood is incompressible. Considering the total volume of the system under consideration,
the quantity of fluids contained in this volume is relatively small, such that the effect of gravity can be neglected.
Since the blood velocity is about $0.1\;\unit{mm/s}$ in the arterioles and venules and about $0.01\;\unit{mm/s}$ in the
capillary bed in a human system \cite[Tab. 1.7]{formaggia2010cardiovascular}, 
it can be concluded that the Reynolds numbers in the whole microvascular network are significantly
lower than $1.0$. Therefore modeling assumption (A4) is reasonable \cite{fung1971microcirculation}. The pulsatility of blood flow can be neglected (A5),
due to the fact that the Womersley numbers in the arterioles, venules and the capillaries are much smaller than $0.1$. The Womersley number is
a dimensionless number relating the frequency of a pulse and the viscosity of a fluid to each other
\cite[Tab. 1.7]{formaggia2010cardiovascular}. The last modeling assumption (A6) is motivated by the fact that the vessel walls of capillaries
consist of a thin layer of endothelial cells with gaps between them such that blood
plasma, oxygen and other substances can migrate into the interstitial
space and then to the tissue cells or vice versa. Contrary to that, the vessel walls of the larger vessels (arterioles and venules) are thicker and not
as permeable, since they exhibit a continuous endothelial cell layer that is surrounded by one or two concentric layers of smooth muscle cells
\cite{despopoulos2003color} \cite[Chap.1]{formaggia2010cardiovascular}.

\section{Models and numerical methods}
\label{sec:Models}

After describing the data set used in this work and discussing the basic modeling assumptions, we outline in this section
the two modeling approaches under consideration, i.e., the fully-discrete 3D-1D model and the hybrid 3D-1D model. In particular, we discuss
coupling concepts between the different parts of the microvascular network and the tissue.
Finally, the numerical discretization of the model equations is shortly explained.

We start by introducing some notation. The cuboid containing the microvascular network depicted
in Fig. \ref{fig:extracted_network} is denoted by $\Omega$ and given by:
\begin{align*}
\Omega = \Big\{\mathbf{x} = \left(x_1,x_2,x_3 \right)^\top \left| \;0\;\unit{mm} < x_1, x_2 < 1.13662\;\unit{mm} \right.\\
\wedge\; 8.75\cdot 10^{-4}\;\unit{mm} < x_3 < 2.16388\;\unit{mm}  \Big\}.
\end{align*}
As it has been described in Section \ref{sec:ProblemSetting}, we assume that the vascular system $\Lambda$ under consideration has 
been segmented and approximated by extracting the
midline of each blood vessel and approximating the midlines by straight segments $\Lambda_k\subset \Lambda$. 
Each segment is equipped with a constant radius value $R_k$ for $k=1,...,N$.
Therefore, the entire vascular system is assumed to be given by the union of $N$ cylinders $V_k$ of radius $R_k$ for $k=1,...,N$. 
We then split the domain in two parts:
$$
\Omega_{\mathrm{v}}=\bigcup_{i=k}^N V_k\qquad\text{and}\qquad \Omega_{\mathrm{t}}=\Omega\setminus \Omega_{\mathrm{v}},
$$
where $\Omega_{\mathrm{v}}$ represents the vascular system and $\Omega_{\mathrm{t}}$ the tissue. 
The entire vascular system can now be described by a graph, whose edges are represented by the center lines 
$\Lambda_k$ of each cylinder $V_k$, for $k=1,...,N$. Each segment $\Lambda_k$ is parametrized by the arc length 
$s_k$, and $\mathbf{\lambda}_k$ is the tangent unit vector
determining the orientation of the centerline of $V_k$. Let us denote the two endpoints of the edge $\Lambda_k$ by $\mathbf{x}_{k,1},\mathbf{x}_{k,2}\in\Omega$.
With this notation at hand, each cylinder $V_k$ can be defined as:
$$
V_k=\left\{\mathbf{x} \in \Omega_{\mathrm{v}} \left| \; \mathbf{x}=\mathbf{x}_{k,1}+ s_k \cdot  \bs{\lambda}_k + \mathbf{r}_k \right. \right\},
$$
where $\mathbf{x}_{k,1} + s_k \cdot \bs{\lambda}_k \in \Lambda_k=\mathcal{M}_k(\Lambda'\subset\R^1)$ and
$$
\mathbf{r}_k \in\mathcal{D}_{\Lambda_k}(R_k)=\left\{r\mathbf{n}_{\Lambda_k}\left( s_k, \theta \right)\;:\; r\in \left[0,R_k\right),
\;s_k \in \left(0,\left| \mathbf{x}_{k,2}-\mathbf{x}_{k,1} \right| \right),\; \theta \in \left[0,2\pi \right) \right\}.
$$
$\mathcal{M}_k$ is a mapping from a reference domain $\Lambda'$ to the manifold $\Lambda_k \subset \Omega$, 
and $\mathbf{n}_{\Lambda_k}$ denotes the set of unit
normal vectors with respect to $\Lambda_k$. $\Gamma_k$ is the lateral surface of the branch $V_k$, and the union $\Gamma$ corresponds approximately to the
surface of the vascular system, since the cylinders $V_k$ may not perfectly match to each other:
$\Gamma = \bigcup_{k=1}^N \Gamma_k.$ The given blood vessel network is separated into two parts:
Based on a fixed threshold $R_{\mathrm{T}}$, we subdivide the network $\Lambda$ in two subsets $\Lambda_{\mathrm{L}}$ consisting of large vessels, 
which are considered to be impermeable and
$\Lambda_{\mathrm{C}}$ consisting of capillaries. Associated to this separation, we define two index sets $I_{\mathrm{L}},\;I_{\mathrm{C}}:$
$$
I_{\mathrm{L}} := \left\{ k \in \left\{ 1,\ldots,N \right\}\;\left|\; R_k\geq R_{\mathrm{T}} \right. \right\} \quad\text{and}\quad
I_{\mathrm{C}}:= \left\{ k \in \left\{ 1,\ldots,N \right\}\;\left|\; R_k<R_{\mathrm{T}} \right. \right\}.
$$
Using these definitions, $\Lambda_{\mathrm{L}}$ and $\Lambda_{\mathrm{C}}$ can be represented as follows:
$$
\Lambda_{\mathrm{L}}:= \bigcup_{k \in I_{\mathrm{L}}} \Lambda_k \quad\text{and}\quad
\Lambda_{\mathrm{C}}:= \bigcup_{k \in I_{\mathrm{C}}} \Lambda_k.
$$
The surface of the capillaries is given by:
$\Gamma_{\mathrm{C}}:= \bigcup_{k \in I_{\mathrm{C}}} \Gamma_k$. Finally, to each node of the graph a pressure value is assigned.
The computation of the pressure values corresponding to the interior nodes is outlined in Section 3.3, while at the boundary nodes the 
pressure values are prescribed from the given data set (see Section 2).

\subsection{Fully-discrete 3D-1D model}
With respect to the surrounding tissue, we assume that it can be considered as a porous medium with 
an isotropic scalar permeability $K_{\mathrm{t}}\;\left[ m^2 \right]$ \cite{khaled2003role}.
In order to model the flow in a porous medium, one can use the Darcy equation to obtain a pressure field \cite{bear2013dynamics,helmig1997multiphase}.

Within the 1D network $\Lambda$, we assume that the flow is described by Hagen-Poiseuille's law with a mass balance equation.
Differentiation over the branches is defined using the tangent unit vector as
$d w/d s_k = \nabla w\cdot \bs{\lambda}_k\,$ on $ \Lambda_k$, i.e. $d /d s_k$ represents the projection of $\nabla$
along $\bs{\lambda}_k$.
The governing flow equations for the whole network $\Lambda$ are obtained by summing the governing equation over the 
index $k$ and introducing the global variable $s$.
Furthermore, we assume that the mass transfer from the vessel into the tissue matrix and vice versa occurs across 
the membrane $\Gamma_{\mathrm{C}}$ of the capillaries according to Starling's filtration law \cite{levick2010microvascular,starling1896absorption}.
Following the approach presented in \cite{koppl2018mathematical}, these exchanges processes are accounted for by a Dirac measure $\delta_{\Gamma_{\mathrm{C}}}$
in the source term of the 3D tissue problem:
\begin{equation}
\label{eq:fully_coupled_system}
\left\{
\begin{aligned}
-\nabla \cdot \left( \rho_{\mathrm{int}} \frac{K_{\mathrm{t}}}{\mu_{\mathrm{int}}} \nabla p^{\mathrm{t}}\right)&=
L_{\mathrm{cap}} \rho_{\mathrm{int}}  \left(p^{\mathrm{v}}-\overline{p^{\mathrm{t}}} - \left( \pi_{\mathrm{p}}
- \pi_{\mathrm{int}} \right) \right)\delta_{\Gamma_{\mathrm{C}}},
&\quad&\text{in } \Omega,\\
\ \\
-\frac{d}{d s}\left( \rho_{\mathrm{bl}} \pi R^2 \frac{K_{\mathrm{v}}}{\mu_{\mathrm{bl}}} \frac{d p^{\mathrm{v}}}{d s} \right)&=
2 \pi R L_{\mathrm{cap}} \rho_{\mathrm{int}} \left(\overline{p^{\mathrm{t}}}-p^{\mathrm{v}} + \left( \pi_{\mathrm{p}} 
- \pi_{\mathrm{int}} \right)\right),&\quad&\text{in } \Lambda,\\
\ \\
\rho_{\mathrm{int}} \frac{K_{\mathrm{t}}}{\mu_{\mathrm{int}}} \nabla p^{\mathrm{t}} \cdot \mathbf{n} &=0,&\quad&\text{on } \partial\Omega,\\
p^{\mathrm{v}}&=p^{\mathrm{v}}_{\mathrm{D}},&\quad&\text{on }\partial\Lambda,
\end{aligned}\right.
\end{equation}
For each edge $\Lambda_k$ and its lateral surface $\Gamma_k$ and $f\in L^2(\Gamma_k)$, we indicate with
$f\delta_{\Gamma_k}$ the linear operator in $\mathcal{C}(\Omega)$ defined by
$$
\langle f\delta_{\Gamma_k},\phi\rangle= \int_{\Lambda_k} \int_{\partial B\left(s_k,R_k\right)} f\phi\;dS\;ds_k,\;\forall \phi\in\mathcal{C}^\infty_0(\Omega),
\; k \in I_{\mathrm{C}}.
$$
The symbol $\partial B(s_k,R_k)$ denotes the circle with center in $\Lambda_k(s_k)$ and perpendicular to $\bs{\lambda}_k$. 
The variable $p^{\mathrm{v}}$ denotes the pressure in the network and $p^{\mathrm{t}}$
the pressure in the tissue. The term $\overline{p^{\mathrm{t}}} (s_k)$ in~\eqref{eq:fully_coupled_system} represents the average value of 
$p^{\mathrm{t}}$ with respect to the circle $\partial B(s_k,R_k)$:
\begin{equation}
\label{eq:average}
\overline{p^{\mathrm{t}}}(s_k)=\frac{1}{2\pi R_k} \int_0^{2\pi}p^{\mathrm{t}} \left( \Lambda_k \left( s_k \right)
+R_k\mathbf{n}_{\Lambda_k}(s_k,\theta) \right)R_k\;d\theta,\qquad k \in I_{\mathrm{C}}.
\end{equation}
$\pi_{\mathrm{p}}$ and $\pi_{\mathrm{int}}$ represent the oncotic pressures of the plasma and the interstitial fluid, respectively. 
According to \cite{levick1991capillary}, these values are assumed to be constant. The hydraulic conductivity 
$L_{\mathrm{cap}}\;\left[ \unit{m/(Pa \cdot s)}\right]$
of the membrane is assumed to be constant for the segments $\Lambda_k,\;k \in I_{\mathrm{C}}$. 
The radius $R\;\left[\unit{m} \right]$ and the permeability $K_{\mathrm{v}}\;\left[\unit{m}^2 \right]$ are defined for
each segment $\Lambda_k$ as:
\begin{equation}
\label{eq:permeability_vessel}
\left. R \right|_{\Lambda_k} = R_k \text{ and } \left. K_{\mathrm{v}} \right|_{\Lambda_k} = \frac{ R_k^2}{8}.
\end{equation}
The viscosity $\mu_{\mathrm{bl}}\;[\unit{Pa\cdot s}]$ of blood is computed according to \eqref{eq:viscosity}, 
while $\mu_{\mathrm{int}}\;[\unit{Pa\cdot s}]$ represents the viscosity of interstitial fluid. $\rho_{\mathrm{int}}\;\left[\unit{kg/m^3} \right]$
and $\rho_{\mathrm{bl}}\;\left[\unit{kg/m^3} \right]$ are the densities of the interstitial fluid and of blood, respectively. 
$p^{\mathrm{v}}_{\mathrm{D}}\;\left[ \unit{Pa}\right]$ are the boundary data,
obtained from the data set previously described in Section~\ref{sec:ProblemSetting}. On the boundary of the tissue continuum 
$\partial \Omega \subseteq \partial \Omega_{\mathrm{t}} \setminus \Gamma_{\mathrm{C}} $,
we set homogeneous Neumann boundary conditions. Due to the fact that we have no boundary data available homogeneous Neumann 
boundary conditions have been chosen. However, if measurements or other data become available,
this boundary condition can easily be replaced.

Alternative 3D-1D PDE-systems simulating flow in microvascular networks, can be found in \cite{d2012finite,d2008coupling}. 
The difference to the presented coupling approach is that in the source term of
the tissue problem, the Dirac measure is concentrated on the midlines of the vascular system. As a result singularities along 
the network midlines are introduced in the 3D pressure field
\cite{koppl2016local,gjerde2018splitting}. Moreover, there are no estimates for the modeling errors arising from this type of coupling concepts. 
In \cite{koppl2018mathematical} we proved for a two-dimensional
model problem that the coupling approach in \eqref{eq:fully_coupled_system} causes a small modeling error, if the radii of the network are small 
compared to the considered tissue matrix. Furthermore, the pressure
field in the 3D tissue matrix does not exhibit any singularities, but only kinks at vessel surfaces.

\subsection{Hybrid 3D-1D model}

\begin{figure}[h!]
\centering
\includegraphics[width=0.92\textwidth]{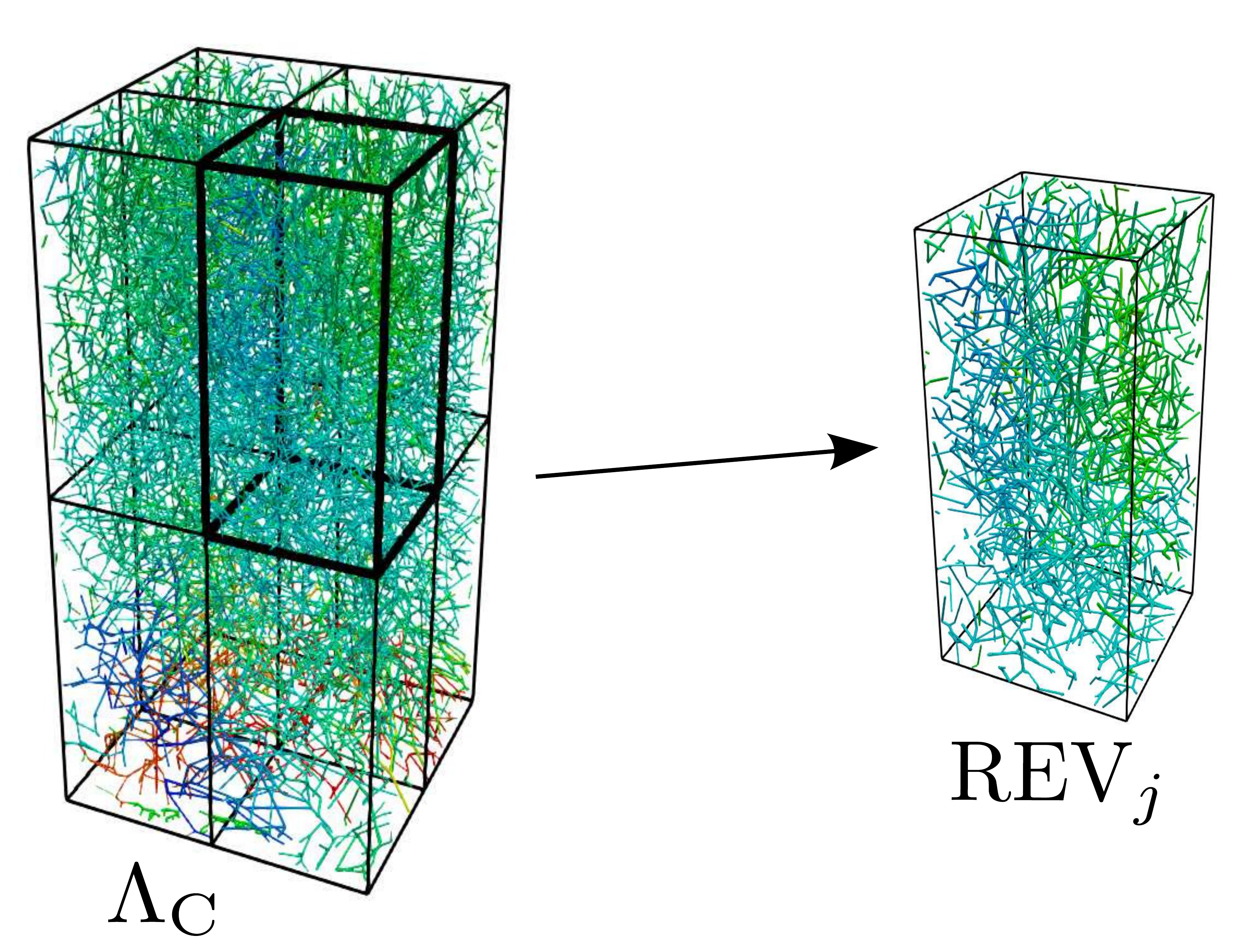}
\caption{\label{fig:capillaries_subdivision} Subdivision of the capillaries $\Lambda_\mathrm{C}$ into eight REVs.}
\end{figure}

\begin{figure}[!t]
\centering
\includegraphics[width=0.7\textwidth]{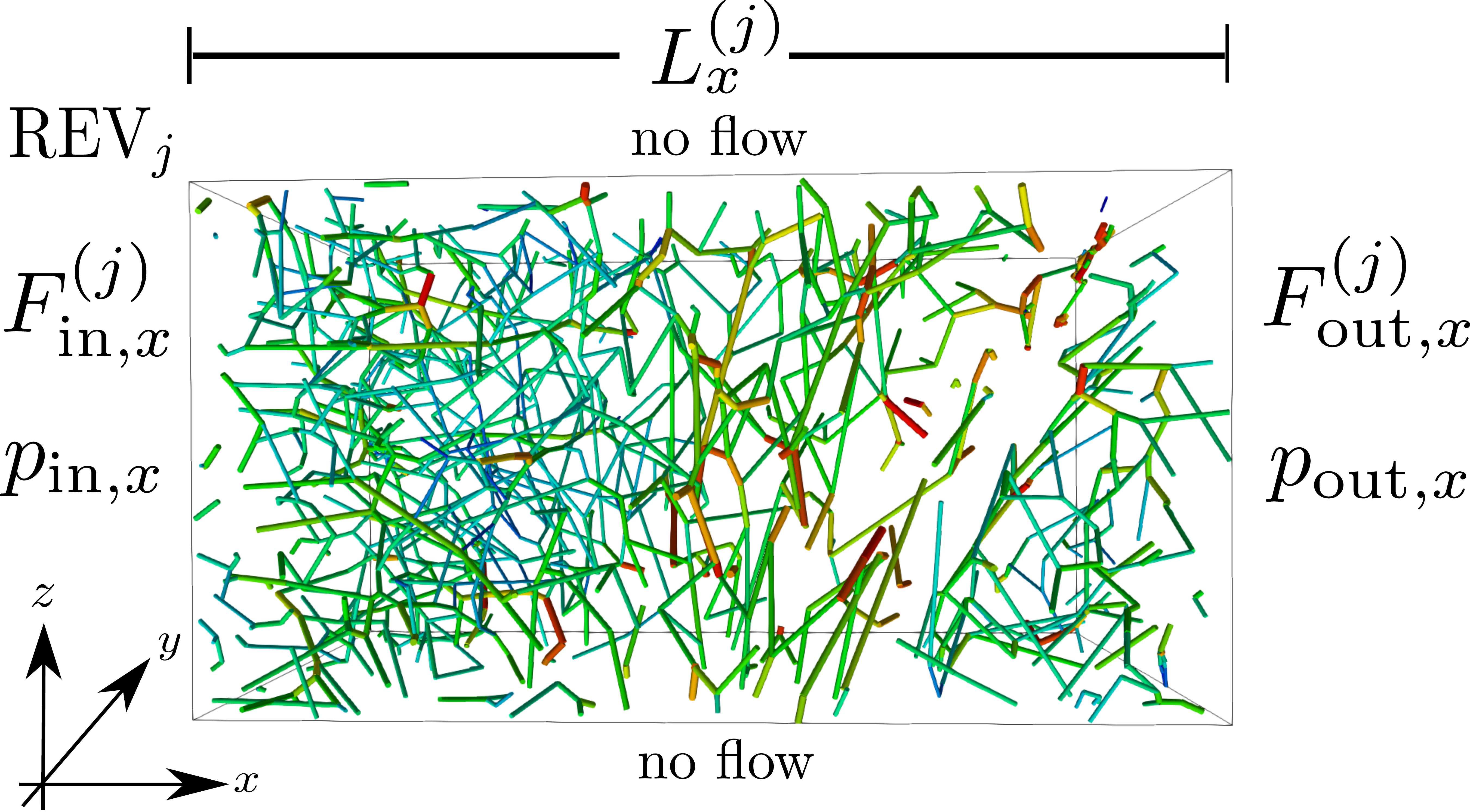}
\caption{\label{fig:permeability_comp} Homogenization of a capillary network contained in a $\text{REV}_j$. 
On the facets $F_{\mathrm{in},x}^{(j)}$ and $F_{\mathrm{out},x}^{(j)}$
whose normals are aligned with the $x$-axis, we apply pressures $p_{\mathrm{in},x}$ and $p_{\mathrm{out},x}$, 
while on the other facets a no-flow condition is imposed. Computing
the flux between $F_{\mathrm{in},x}^{(j)}$ and $F_{\mathrm{out},x}^{(j)}$, the permeability component $k_x^{(j)}$ can be estimated.}
\end{figure}

\begin{figure}[h!]
\begin{center}
\includegraphics[width=0.9\textwidth]{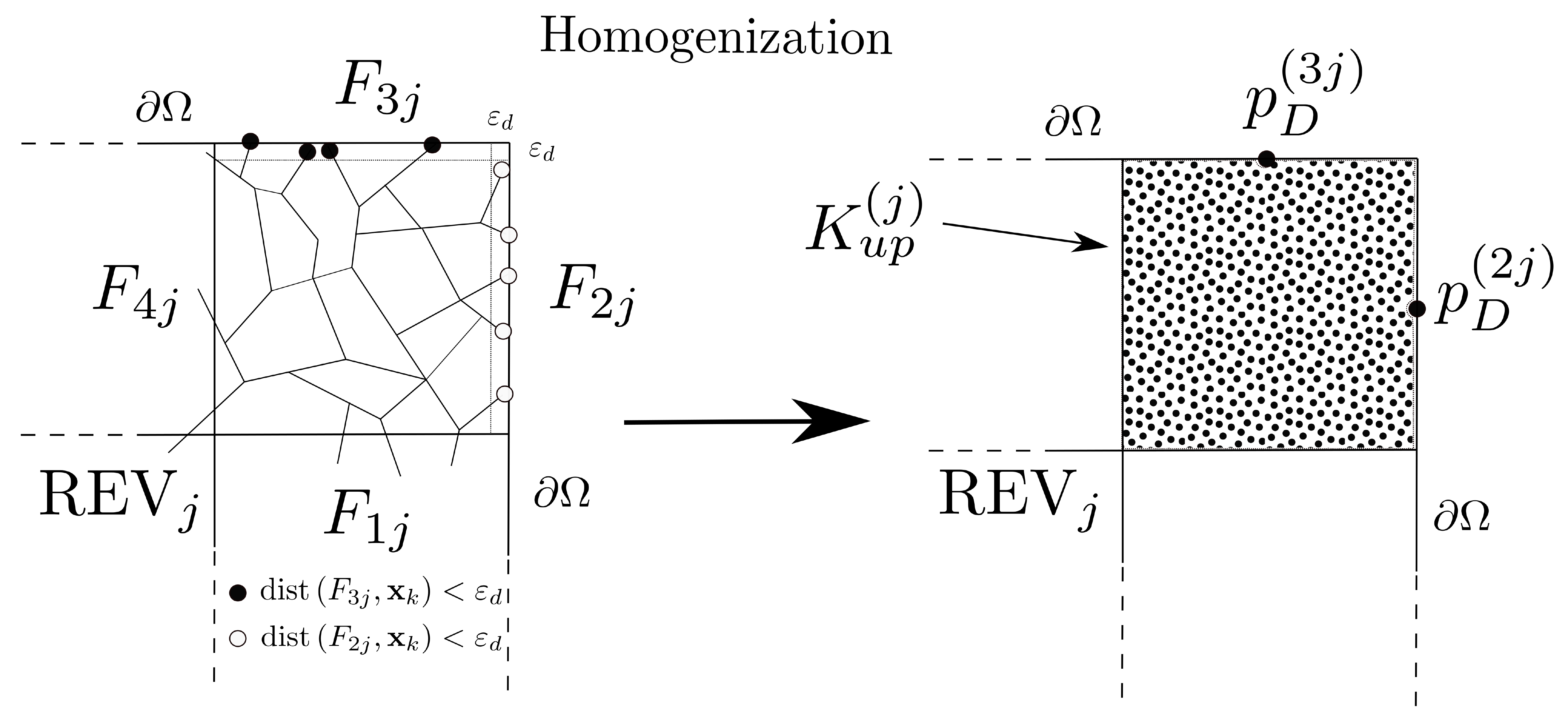}
\end{center}
\caption{\label{fig:boundary_conditions} The figure shows a two dimensional layout parallel to the $x$-$y$ plane
of an $\mathrm{REV}_j$ that is located at an edge of the domain $\Omega$. At the top of the figure,
the discrete capillary network contained in $\mathrm{REV}_j$ is shown,
while in the bottom part of the figure, the homogenized system with the tensor $K_{\mathrm{up}}^{(j)}$ and the averaged boundary pressures
$p_D^{(2j)}$ and $p_D^{(3j)}$ can be seen.}
\end{figure}

Considering the microvascular network in Fig.~\ref{fig:extracted_network}, it can be seen that it consists of
venules and arterioles and a large number of small capillary vessels forming a dense structure.
Due to that, we consider the capillaries as a 3D porous medium, while the venules and arterioles are still considered as
1D vessel systems. This approach has the clear advantage that it does not require
a high-resolution description of the microvascular network under consideration and that it preserves the hierarchy of the larger vessels.
All in all, there are now two coupled 3D continua in $\Omega$, one for the capillary bed and another
one for the tissue. This means that the
hybrid 3D-1D model is a \emph{double-3D continuum} model contrary to the \emph{fully-discrete} model, which is a \emph{single-3D continuum} model.
As for the tissue, the flow in the homogenized capillaries can be described using Darcy's law, and
we pose the following problem for the corresponding unknown $p^{\mathrm{cap}}$:
\begin{equation}
\label{eq:darcy_capillaries}
\left\{
\begin{aligned}
-\nabla \cdot \left( \rho_{\mathrm{bl}} \frac{K_{\mathrm{up}}}{\mu_{\mathrm{bl}}^{\mathrm{up}}} \nabla p^{\mathrm{cap}}\right)&= 
q^{\mathrm{cap}},&\quad& \text{in } \Omega,\\
p^{\mathrm{cap}}&=p^{\mathrm{cap}}_{\mathrm{D}},&\qquad&\text{on } \partial\Omega.
\end{aligned}\right.
\end{equation}
$K_{\mathrm{up}}$ is the corresponding permeability tensor, $\mu_{\mathrm{bl}}^{\mathrm{up}}$ is an averaged viscosity, $q^{\mathrm{cap}}$ 
indicates the source term
and $p^{\mathrm{cap}}_{\mathrm{D}}$ denotes the Dirichlet boundary condition. These terms and parameters are described in the following.
\\ \\
\textbf{Computing the tensor} $K_{\mathrm{up}}$ \textbf{and the averaged viscosity} $\mu_{\mathrm{bl}}^{\mathrm{up}}$
\\ \\
Let us assume that the domain $\Omega$ can be decomposed into \emph{representative elementary volumes} (REVs) \cite{helmig1997multiphase}, that is:
\begin{equation}\label{eq:omega_rev_sub}
\overline{\Omega}=\bigcup_{j=1}^{N_{\mathrm{REV}}} \overline{\mathrm{REV}_j},
\end{equation}
where $N_{\mathrm{REV}}$ is the total number of REVs. Furthermore, we assume that each $\mathrm{REV}_j \subset \Omega$ 
is a rectangular cuboid, as depicted in Fig. \ref{fig:capillaries_subdivision}. With respect to $\mathrm{REV}_j$, the viscosity 
$\mu_{\mathrm{bl}}^{\mathrm{up}}$ is defined as follows:
$$
\mu_{\mathrm{bl}}^{\mathrm{up}}(\mathbf{x}) = \mu_{\mathrm{bl},j}^{\mathrm{up}}, \text{ if } \mathbf{x} \in \mathrm{REV}_j, \;\text{ where }\;
\mu_{\mathrm{bl},j}^{\mathrm{up}} = \frac{1}{\left| I_{\mathrm{C},j} \right|} \sum_{k \in I_{\mathrm{C},j}} \mu_{\mathrm{bl}} \left(2 \cdot R_k / \mu m \right).
$$
The viscosity $\mu_{\mathrm{bl}}$ is given by \eqref{eq:viscosity}, and the set $I_{\mathrm{C},j}$ is defined as follows:
$$
I_{\mathrm{C},j} = \left\{ k \in I_{\mathrm{C}} \left|\; \Lambda_k \cap \mathrm{REV}_j \neq \emptyset \right. \right\}.
$$
Now, we determine the tensor $K_{\mathrm{up}}$ which represents the permeability of the homogenized capillary bed.
In~\cite{el2015multi}, the authors computed full permeability tensors for periodic cells based on computer-generated capillary networks,
whose properties statistically match measurements in brain tissue. 
The off-diagonal entries were found to be on average two orders of magnitude smaller than the diagonal terms~\cite[Table 2]{el2015multi}.
Following this observation, we make the following simplifying assumptions for $K_{\mathrm{up}}$:\\
\begin{itemize}
 \item[(i)] the permeability tensor $K_{\mathrm{up}}$ is constant on each $\mathrm{REV}_j$:
 $$
 K_{\mathrm{up}}(\mathbf{x})= K_{\mathrm{up}}^{(j)}, \text{ if } \mathbf{x} \in \mathrm{REV}_j.
 $$
 \item[(ii)] the permeability tensor $K_{\mathrm{up}}^{(j)}$ is diagonal:
 $$
 K_{\mathrm{up}}^{(j)} =\begin{pmatrix}
 k_x^{(j)} & 0         & 0 \\
 0        & k_y^{(j)} & 0 \\
 0        &         0 & k_z^{(j)}
 \end{pmatrix},
 $$
with $k_x^{(j)},k_y^{(j)},k_z^{(j)}>0$.
\end{itemize}
In order to determine the components of the permeability tensor $K_{\mathrm{up}}^{(j)}$, we apply the upscaling strategy presented in \cite{reichold2009vascular}
which we briefly describe in the following. For simplicity, we restrict ourself to the computation of the permeability component $k_x^{(j)}$. The other
permeability components $k_y^{(j)}$ and $k_z^{(j)}$ can be computed in an analogous way. As a first step to compute
this quantity, we apply a no-flow boundary condition to all the facets of the $\text{REV}_j$ whose face normals are not aligned with the $x$-axis
(see Fig. \ref{fig:permeability_comp}).
The remaining facets are denoted by $F_{\mathrm{in},x}^{(j)}$ and $F_{\mathrm{out},x}^{(j)}$. Between these facets a pressure gradient 
is applied by imposing a pressure
$p_{\mathrm{in},x}$ on $F_{\mathrm{in},x}^{(j)}$
and a pressure $p_{\mathrm{out},x}$ on $F_{\mathrm{out},x}^{(j)}$, where $p_{\mathrm{in},x} > p_{\mathrm{out},x}$. 
This results into a volume flux from $F_{\mathrm{in},x}^{(j)}$ to $F_{\mathrm{out},x}^{(j)}$.
Using the Vascular Graph Model (VGM) described in Subsection \ref{sec:numerical_discretization}, we compute the pressure field in
$$
\Lambda_{\mathrm{C},j} = \bigcup_{k \in I_{\mathrm{C},j}} \Lambda_k.
$$
By means of this pressure field, the volume flux $VF_{\mathrm{out},x}^{(j)}$ through $F_{\mathrm{out},x}^{(i)}$ is computed as follows:
$$
VF_{\mathrm{out},x}^{(j)}=\sum_{\mathbf{x}_k \in \Lambda_{\mathrm{C},j} \cap F_{\mathrm{out},x}^{(j)}}
\frac{\pi R_k^2 K_{\mathrm{v}}\left( s_k \right)}{\mu_{\mathrm{bl}}}\cdot
\frac{\partial p^{\mathrm{v}}}{\partial s}\left( s_k \right),
$$
where $\mathbf{x}_k = \Lambda_k\left(s_k \right)$ for $k \in I_{\mathrm{C},j}$. Based on $VF_{\mathrm{out},x}^{(j)}$, 
the permeability component $k_x^{(j)}$ is approximated as follows:
\begin{equation}\label{eq:k_up}
k_x^{(j)} \approx \frac{VF_{\mathrm{out},x}^{(j)} \cdot \mu_{\mathrm{bl},j}^{\mathrm{up}} 
\cdot L_x^{(j)}}{L_y^{(j)} \cdot L_z^{(j)} \cdot \left(p_{\mathrm{in},x}-  p_{\mathrm{out},x}\right)}.
\end{equation}
$L_x^{(j)},\;L_y^{(j)}$ and $L_z^{(j)}$ are the edge lengths of $\mathrm{REV}_j$ in
the $x$-, $y$- and $z$-direction, respectively.
An open issue in this context is the choice of the $\mathrm{REV}$-size. In order to clarify this issue,
we refer to Subsection~\ref{sec:K_up}, in which the admissible size of the $\mathrm{REV}$s is determined numerically.
\\ \\
\textbf{Computing the boundary conditions} $p^{\mathrm{cap}}_{\mathrm{D}}$
\\ \\
For the computation of the boundary data $p^{\mathrm{cap}}_{\mathrm{D}}$, we consider the $\mathrm{REV}_j$ that are adjacent 
to the boundary $\partial \Omega$, i.e., $\overline{\mathrm{REV}}_j \cap \partial \Omega \neq \emptyset$. The six facets of a $\mathrm{REV}_j$ 
are denoted by $F_{1j}, \ldots, F_{6j}$.
For each $F_{ij}$ with $F_{ij} \subset \partial \Omega$, we compute an averaged boundary value $p_{\mathrm{D}}^{(ij)}$. 
This is done, by averaging all the pressure values $p_k$
that are assigned to the boundary nodes $\mathbf{x}_k \in \partial \Lambda_{\mathrm{C}}$ and whose distance to the face $F_{ij}$ 
is smaller than a small parameter $\varepsilon_d>0$ (see Fig.~\ref{fig:boundary_conditions}): 
$\text{dist} \left(F_{ij},\mathbf{x}_k \right) < \varepsilon_d$. For the rest of the paper, we set $\varepsilon_d=10^{-8}\;\unit{m}$. 
Assuming that $N_{ij}$ nodes are fulfilling these conditions, we compute the averaged pressure by an arithmetic mean:
$$
p_{\mathrm{D}}^{(ij)} = \frac{1}{N_{ij}} \sum_{k=1}^{N_{ij}} p_k.
$$
In order to obtain a smooth function for the boundary values $p^{\mathrm{cap}}_{\mathrm{D}}$, we assign the pressure values $p_{\mathrm{D}}^{(ij)}$
to the centers of the faces $F_{ij}$ and construct a linear interpolant $p^{\mathrm{cap}}_{\mathrm{D}}$ on $\partial \Omega$ based on the described 
setup by means of the function \texttt{interpolate}
from the package PDELab of Dune~\cite{blatt2016distributed}.
\\ \\
\textbf{Computing the source term} $q^{\mathrm{cap}}$
\\ \\
The source term $q^{\mathrm{cap}}$ in \eqref{eq:darcy_capillaries} can be splitted into a component $q^{\mathrm{cap}}_{\mathrm{v}}$
accounting for the impact of the larger vessels and a component $q^{\mathrm{cap}}_{\mathrm{t}}$ incorporating the influence of the surrounding tissue:
\begin{equation}
\label{eq:source_cap}
q^{\mathrm{cap}} = q^{\mathrm{cap}}_{\mathrm{v}} + q^{\mathrm{cap}}_{\mathrm{t}}.
\end{equation}
Let $\mathbf{x}_k$ be a node at the boundary of the subset $\Lambda_{\mathrm{L}}$ and in the interior of the domain $\Omega$,
i.e., $\mathbf{x}_k \in \partial \Lambda_{\mathrm{L}} \cap \Omega$. Furthermore, we assume that $\mathbf{x}_k \in \Lambda_k\cap \mathrm{REV}_j$ for an
index $k\in I_{\mathrm{L}}$ and an index $j\in\{1,...,N_{\mathrm{REV}}\}$. This means, there exists a $s_k$ such that $\mathbf{x}_k=\Lambda_k(s_k)$.
Then the flux occurring between the edge $\Lambda_k$ and the capillary continuum is given by:
\begin{equation}
\label{eq:BCLargeVessels}
\rho_{\mathrm{bl}} \frac{\pi R_k^2 K_{\mathrm{v}}\left( s_k \right)}{\mu_{\mathrm{bl}}}\cdot \frac{\partial p^{\mathrm{v}}}{\partial s}\left( s_k \right) =
\alpha_{\mathrm{v}}^{\mathrm{cap}}\left( \Lambda_k,\mathrm{REV}_j \right) \left(  p^{\mathrm{cap}}_{(j)} -p^{\mathrm{v}}\left( s_k \right) \right),
\end{equation}
where the factor $\alpha_{\mathrm{v}}^{\mathrm{cap}}$ and the averaged pressure $p^{\mathrm{cap}}_{(j)}$ are defined as follows:
\begin{equation}
\label{eq:REV_pressure_cap}
\alpha_{\mathrm{v}}^{\mathrm{cap}}\left( \Lambda_k,\mathrm{REV}_j \right)
= \rho_{\mathrm{bl}}\cdot\frac{ \pi R_k^2 K_{\mathrm{v}}^{(j)}}{\mu_{\mathrm{bl}}^{\mathrm{up}}\ell_{\mathrm{c}}^{(kj)}}
\;\text{ and }\; p^{\mathrm{cap}}_{(j)} = \frac{1}{\left| \mathrm{REV}_j \right| } \int_{\mathrm{REV}_j} p^{\mathrm{cap}}\left(x \right)\;dx.
\end{equation}
The coefficient $K_{\mathrm{v}}^{(j)}$ represents the permeability of the capillaries connected to the coupling point $\mathbf{x}_k$.
The parameter $\ell_{\mathrm{c}}^{(kj)}\;\left[ \unit{m} \right]$ indicates the average length of the blood flow paths
that begin at $\mathbf{x}_k$ and are contained in $\mathrm{REV}_j$ which is not known a priori and has to be estimated. Therefore, we set:
\begin{equation}\label{eq:lc_estimate}
\frac{K_{\mathrm{v}}^{(j)}}{\ell_{\mathrm{c}}^{(kj)}}=\alpha \frac{\overline{K_{\mathrm{v}}^{(j)}}}{L_j},
\end{equation}
where $L_j$ is the smallest edge length of the $\mathrm{REV}_j$, and $\overline{K_{\mathrm{v}}^{(j)}}$ denotes the arithmetic average of the permeabilities
\eqref{eq:permeability_vessel} of the capillaries contained in the $\mathrm{REV}_j$.  A numerical study to determine the optimal 
value of the parameter $\alpha\in (0,1)$
(with respect to the fluxes within the system) for the problem under consideration is postponed to Subsection \ref{sec:mf_comp}.
Considering the right hand side of \eqref{eq:BCLargeVessels} one has to note that the term
$$
\frac{p^{\mathrm{cap}}_{(j)} -p^{\mathrm{v}}\left( s_k \right)}{\ell_{\mathrm{c}}^{(kj)}}
$$
represents a finite difference approximation of a pressure gradient at an outlet of the larger vessels. According to the REV-concept in porous media theory
\cite{bear2013dynamics,helmig1997multiphase}, one can assign to each REV an averaged pressure $p^{\mathrm{cap}}_{(j)}$, which stands for the pressure field
in the $\mathrm{REV}_j$. The finite difference above can be considered as
an approximation of the pressure gradient between an outlet of $\Lambda_{\mathrm{L}}$ and the homogenized capillary bed in $\mathrm{REV}_j$.
Using \eqref{eq:BCLargeVessels}, the source term $q^{\mathrm{cap}}_{\mathrm{v}}$ is computed for $\mathbf{x} \in \mathrm{REV}_j$ by:
\begin{equation}
\label{eq:source_cap_l}
q^{\mathrm{cap}}_{\mathrm{v}} \left( \mathbf{x} \right) =
\sum_{\mathbf{x}_k \in \partial \Lambda_{\mathrm{L}} \cap \mathrm{REV}_j} 
\frac{\alpha_{\mathrm{v}}^{\mathrm{cap}}\left( \Lambda_k,\mathrm{REV}_j\right)}{\left| \mathrm{REV}_j \right|}
\left( p^{\mathrm{v}}\left( s_k \right) - p^{\mathrm{cap}}_{(j)} \right),
\end{equation}
such that the model is mass conservative.
\ \\ \\
\textbf{Computing the source term} $q_{\mathrm{t}}^{\mathrm{cap}}$
\ \\ \\
As in the case of the fully-discrete 3D-1D model, the tissue is considered as a porous structure.
The main difference to the capillaries is that we assume an isotropic permeability $K_{\mathrm{t}}$ for the tissue.
Furthermore the interstitial fluid is assumed to have
a constant viscosity $\mu_{\mathrm{int}}$. By means of Darcy's law, the pressure $p^{\mathrm{t}}$ can be computed by:
\begin{equation}
\label{eq:tissue}
\left\{
\begin{aligned}
 -\nabla \cdot \left( \rho_{\mathrm{int}} \frac{K_{\mathrm{t}}}{\mu_{\mathrm{int}}} \nabla p^{\mathrm{t}} \right) = 
 -q_{\mathrm{t}}^{\mathrm{cap}}, &\qquad &\text{in } \Omega,\\
 \rho_{\mathrm{int}} \frac{K_{\mathrm{t}}}{\mu_{\mathrm{int}}} \nabla p^{\mathrm{t}} \cdot \mathbf{n} = 0, &\qquad &\text{on } \partial \Omega.
\end{aligned}\right.
\end{equation}
It remains to specify the source term $q_{\mathrm{t}}^{\mathrm{cap}}$ modeling the fluid transfer between the capillary bed and the intracellular space.
For this purpose, we use as in \eqref{eq:fully_coupled_system} Starling's filtration law with respect to the
vessel surface area $S_j$ that is contained in an $\mathrm{REV}_j$:
$$
S_j = \sum_{k \in I_{\mathrm{C},j} } \left| \widetilde{\Lambda}_k \right| \cdot 2 \pi R_k,
$$
where $\widetilde{\Lambda}_k\subseteq\Lambda_k$ such that $\Lambda_k\cap \mathrm{REV}_j=\widetilde{\Lambda}_k$. 
Using this parameter, $q_{\mathrm{t}}^{\mathrm{cap}}$ is given by:
\begin{equation}
\label{eq:source_cap_t}
q_{\mathrm{t}}^{\mathrm{cap}} \left( \mathbf{x} \right) = \frac{\rho_{\mathrm{int}} \cdot S_j \cdot L_{\mathrm{cap}}}{\left| \mathrm{REV}_j \right|}
\left( p^{\mathrm{t}}\left( \mathbf{x} \right) - p^{\mathrm{cap}}\left( \mathbf{x} \right) + \left( \pi_{\mathrm{p}} - \pi_{\mathrm{int}} \right) \right),
\; \mathbf{x} \in \mathrm{REV}_j.
\end{equation}
\\
\textbf{Summary of the equations governing the hybrid model}
\\ \\
Summarizing all the previous considerations, the hybrid (double-continuum)
3D-3D-1D model is governed by the following equations:
\ \\
\begin{itemize}
 \item Large Vessels (1D discrete network):
 \begin{equation}
 \label{eq:large_vessels}
 \left\{
 \begin{aligned}
 -\frac{\partial}{\partial s}\left( \rho_{\mathrm{bl}} \cdot \pi R^2 \frac{K_{\mathrm{v}}}{\mu_{\mathrm{bl}}} 
 \frac{\partial p^{\mathrm{v}}}{\partial s} \right)&=0,
 &\qquad&\text{in }\Lambda_{\mathrm{L}}, \\
 p^{\mathrm{v}}&=p^{\mathrm{v}}_{\mathrm{D}},&\qquad&\text{on }\partial\Lambda_{\mathrm{L}}\cap \partial\Omega, \\
 \text{flux term in } & \eqref{eq:BCLargeVessels}, &\qquad&\text{on }\partial \Lambda_{\mathrm{L}}\cap\Omega. \\
 \end{aligned}\right.
 \end{equation}
 \item Capillary bed (3D porous medium):
 \begin{equation}
 \label{eq:capillary_system}
 \left\{
 \begin{aligned}
 -\nabla \cdot \left( \rho_{\mathrm{bl}} \frac{K_{\mathrm{up}}}{\mu_{\mathrm{bl}}^{\mathrm{up}}} \nabla p^{\mathrm{cap}}\right)&= q_{\mathrm{v}}^{\mathrm{cap}}
 + q_{\mathrm{t}}^{\mathrm{cap}},&\qquad& \text{in } \Omega,\\
 p^{\mathrm{cap}}&=p^{\mathrm{cap}}_{\mathrm{D}},&\qquad&\text{on } \partial\Omega.
 \end{aligned}\right.
 \end{equation}
 \item Tissue (3D porous medium):
 \begin{equation}
 \label{eq:tissue_system}
 \left\{
 \begin{aligned}
 -\nabla \cdot \left( \rho_{\mathrm{int}} \frac{K_{\mathrm{t}}}{\mu_{\mathrm{int}}} \nabla p^{\mathrm{t}} \right) &
 = -q_{\mathrm{t}}^{\mathrm{cap}}, &\qquad &\text{in } \Omega,\\
 \rho_{\mathrm{int}} \frac{K_{\mathrm{t}}}{\mu_{\mathrm{int}}} \nabla p^{\mathrm{t}} \cdot \mathbf{n} &= 0, &\qquad &\text{on } \partial \Omega.
 \end{aligned}\right.
 \end{equation}
\end{itemize}
The coupling term $q_{\mathrm{v}}^{\mathrm{cap}}$ is given by~\eqref{eq:source_cap_l}, and describes the interactions between 
the extracted network $\Lambda_{\mathrm{L}}$ and the homogenized capillaries. The other coupling term $q_{\mathrm{t}}^{\mathrm{cap}}$ 
is defined by~\eqref{eq:source_cap_t}, and stands for the exchange between tissue and homogenized capillaries.

\subsection{Numerical discretizations}
\label{sec:numerical_discretization}
Next, we briefly describe a numerical scheme that is used to solve the model equations \eqref{eq:fully_coupled_system} and
\eqref{eq:large_vessels}-\eqref{eq:tissue_system}.
The elliptic PDEs governing the flow within the tissue or the upscaled capillaries
are solved numerically by means of a standard cell-centered finite volume method \cite{helmig1997multiphase}, where the numerical fluxes across
the cell surfaces are approximated by the two-point flux method. The choice of this method is motivated by its
intrinsic local mass conservation, and by the fact that we can work with uniform hexahedral meshes.

\begin{figure}[!b]
\centering
\includegraphics[width=0.65\textwidth]{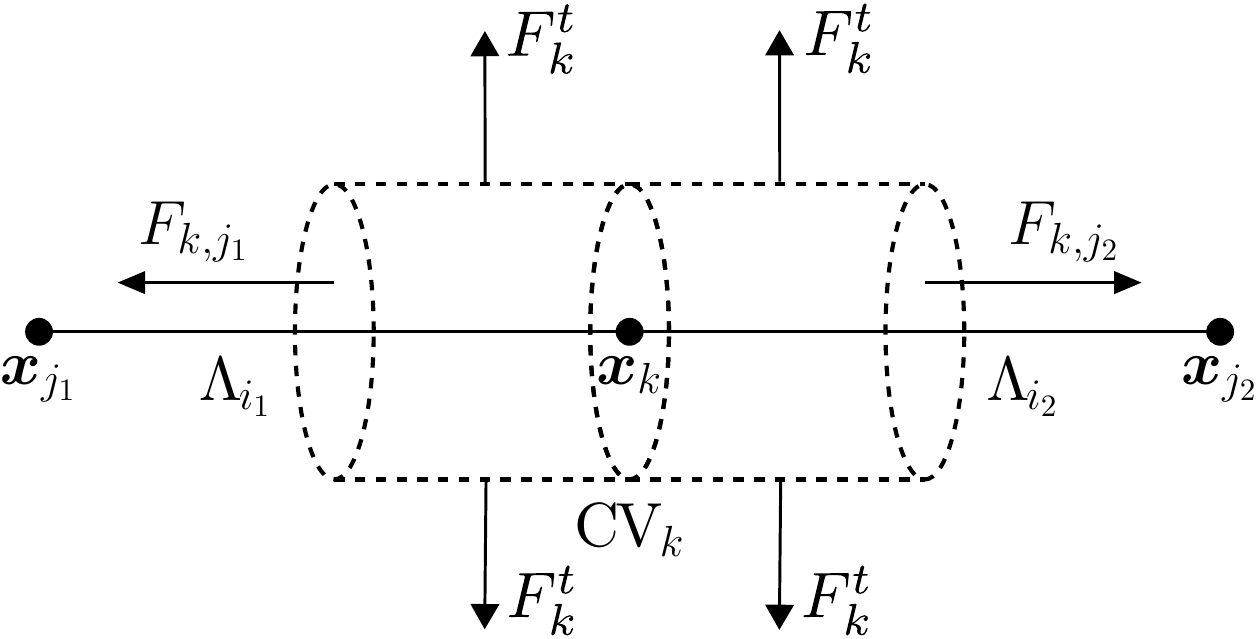}
\caption{Representation of the numerical fluxes through the surface of the
control volume $\mathrm{CV}_k$ with center in $\mathbf{x}_k$. The point shares the edges
$\Lambda_{i_1}$ and $\Lambda_{i_2}$ with the nodes $\mathbf{x}_{j_1}$ and $\mathbf{x}_{j_2}$ of the network, respectively.
\label{fig:control_volume}}
\end{figure}

For the numerical solutions of the network equations in both modeling approaches \eqref{eq:fully_coupled_system} and \eqref{eq:large_vessels},
we employ the vascular
graph model (VGM) \cite{erbertseder2012coupled,peyrounette2018multiscale,reichold2009vascular}. Thereby, we approximate the pressure values at
the grid
nodes $\mathbf{x}_k$ discretizing the network structures. Around each grid node a control volume $\mathrm{CV}_k$ is placed such that the grid
node is in the
center of the control volume (see Fig.~\ref{fig:control_volume}).
In the next step, the fluxes $F_{kj}$ and $F_k^t$
across the surfaces of the control volume are computed and summed up such that the sum of the fluxes is equal to zero:
\begin{equation}
\label{eq:vgm_method}
\sum_{j\in \mathcal{N}(\mathbf{x}_k)} F_{kj}-F_k^t=0,\qquad F_{kj}=\frac{\rho_{\mathrm{bl}} \cdot \pi R_k^4}{8\mu_{\mathrm{bl}}
\left|\mathbf{x}_j-\mathbf{x}_k \right|}\cdot
\Big(  p^{\mathrm{v}}\left(\mathbf{x}_j\right) - p^{\mathrm{v}}\left(\mathbf{x}_k\right) \Big),
\end{equation}
where $R_k$ is the radius of the edge linking $\mathbf{x}_k$ and $\mathbf{x}_j$.
$\mathcal{N}(\mathbf{x}_k)$ denotes the set of indices that share the edge $\Lambda_k$ with the point $\mathbf{x}_k$.

Solving the fully discrete model \eqref{eq:fully_coupled_system}, the walls of the smaller vessels (capillaries) are permeable, 
and a flux across the vessel walls has to be considered, too. This is done by computing an averaged pressure in the tissue with respect 
to the part of the vessel wall touching the control volume, see Fig.~\ref{fig:control_volume}. Then, this pressure value is compared to the network 
pressure associated with the grid node $\mathbf{x}_k$ to determine the flux across the vessel wall:
\begin{align*}
F_k^{\mathrm{t}}=  2 \pi L_{\mathrm{cap}} \cdot \rho_{\mathrm{int}} \left[ \sum_{j \in \mathcal{E}(\mathbf{x}_k) } R_{j}
\int_{CV_k \cap \Lambda_{j}} I \left( p^{\mathrm{v}} \right) - \overline{p^{\mathrm{t}}} 
- \left( \pi_{\mathrm{p}}
- \pi_{\mathrm{int}} \right)\;dS \right].
\end{align*}
$\mathcal{E}(\mathbf{x}_k)$ is the index set for the edges containing the grid node $\mathbf{x}_k$.
$I \left( p^{\mathrm{v}} \right)$ is a linear interpolant for the pressure field $p^{\mathrm{v}}$, whose shape on each edge is 
determined by the two pressure values of the edge.
By this, we obtain for each grid node a mass balance equation. Summarizing these equations, we obtain
a system of equations for the pressure values at the grid nodes.

The reason why we consider a different discretization for the network is the treatment of bifurcations or junctions in the network. Within the VGM
approach, a grid node is placed directly at a bifurcation and the mass balance equation~\eqref{eq:vgm_method} can easily be established. 
Applying the cell-centered
finite volume method, one cannot place a degree of freedom directly at a bifurcation, and therefore the fluxes through the bifurcation point have to be computed.
Since the radii of the branches and the main vessel may be different, the computation of the fluxes requires a careful computation of the numerical
transmissiblity coefficients in~\eqref{eq:vgm_method}. Finally, the continuity of the pressure at a bifurcation is guaranteed using the VGM.

All in all the numerical treatment of the PDEs together with the boundary conditions, source terms and coupling conditions, yields for each model a sparse
linear system of equations. Each block of the system matrix is the discrete representative of an elliptic differential operator or a coupling term,
whereas the contributions of the oncotic pressures and the boundary conditions are incorporated into the right hand side of the system of equations.
For the numerical solution of the linear equation system
a block AMG-preconditioner is applied. The preconditioned system is then solved by a stabilized bi-conjugate gradient method.
This solver was realized using the generic interface of the ISTL-library of DUNE and its AMG implementation \cite{blatt2006iterative,blatt2008generic}.

\section{Numerical tests}
\label{sec:numerical_test}
In this section, we test the numerical models presented in Section \ref{sec:Models} using the data set described in Section~\ref{sec:ProblemSetting}.
Thereby, by means of the fully-discrete 3D-1D model~\eqref{eq:fully_coupled_system} a reference solution is computed 
for the hybrid model~\eqref{eq:large_vessels}-\eqref{eq:tissue_system}.
The results obtained by both approaches are compared with respect to mass 
fluxes and averaged REV-pressures. The flow is driven by the boundary conditions, as described in Sections 3.1. and 3.2. In particular,
no-flow boundary conditions are posed for the flow in the tissue, while prescribed Dirichlet values are considered for the network (see Section 2).
In Table~\ref{table:parameters},
the model parameters for the numerical simulations are summarized. The motivation for the choice of $R_\mathrm{T}=7\;\unit{\mu m}$ is given in 
Subsection~\ref{sec:R_T}, while
the number $N_{\mathrm{REV}}$ in~\eqref{eq:omega_rev_sub} for the hybrid model is discussed in Subsection~\ref{sec:K_up}. Subsections~\ref{sec:mf_comp} and 
\ref{sec:REV_pressures} contain
the numerical results for the mass fluxes and the averaged $\mathrm{REV}$-pressures, respectively. Finally, in Subsection~\ref{sec:boundary}, 
the influence of different boundary conditions are discussed.

\begin{table}[!h]
\centering
\caption{\label{table:parameters} Values of the parameters used for the numerical experiments.}
\begin{center}
\begin{tabular}{|c|c|c|}
\hline
\hline
Discharge hematocrit & $H$ & $0.45$ \\
\hline
Tissue permeability & $K_{\mathrm{t}}$ & $10^{-18}\;\unit{m^2}$\\
\hline
Interstitial fluid viscosity & $\mu_{\mathrm{int}}$ & $1.3 \cdot 10^{-3}\;\unit{Pa\cdot s}$\\
\hline
Plasma viscosity & $\mu_{\mathrm{p}}$ & $1.0 \cdot 10^{-3}\;\unit{Pa\cdot s}$\\
\hline
Blood density & $\rho_{\mathrm{bl}}$ & $1030\;\unit{kg/m^3}$\\
\hline
Interstitial fluid density & $\rho_{\mathrm{int}}$ & $1000\;\unit{kg/m^3}$\\
\hline
Plasma oncotic pressure & $\pi_{\mathrm{p}}$ & $3300\;\unit{Pa}$\\
\hline
Interstitial oncotic pressure & $\pi_{\mathrm{int}}$ & $666\;\unit{Pa}$\\
\hline
Capillary wall hydraulic conductivity & $L_{\mathrm{cap}}$ & $10^{-12}\;\unit{m/(Pa\cdot s)}$\\
\hline
Threshold large vessels/capillaries & $R_{\mathrm{T}}$ & $7\;\unit{\mu m}$\\
\hline
\hline
\end{tabular}
\end{center}
\end{table}

\subsection{Justification of the threshold $R_{\mathrm{T}}$}
\label{sec:R_T}
\begin{figure}[!h]
\centering
\includegraphics[width=0.58\textwidth]{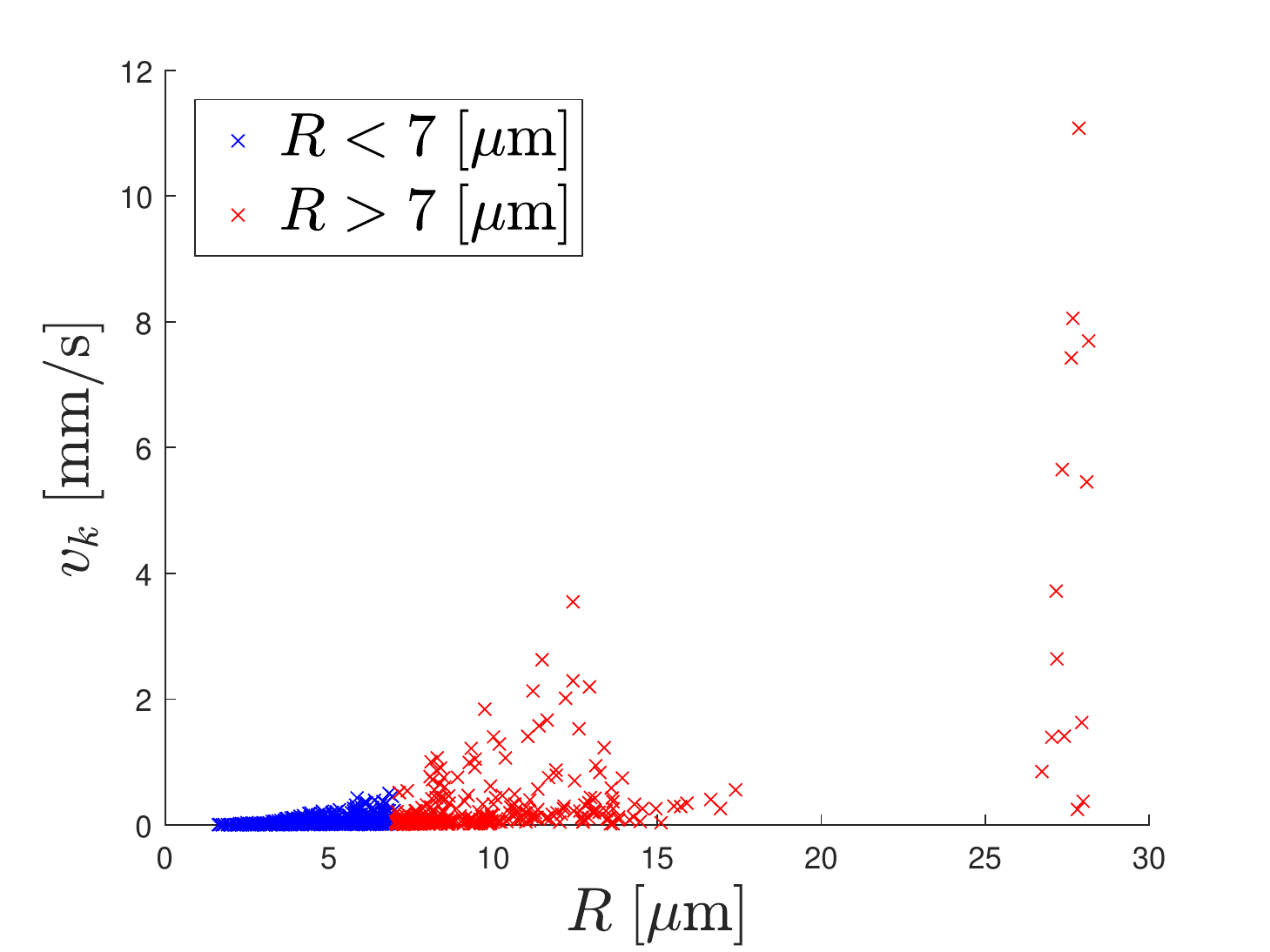}
\includegraphics[width=0.35\textwidth]{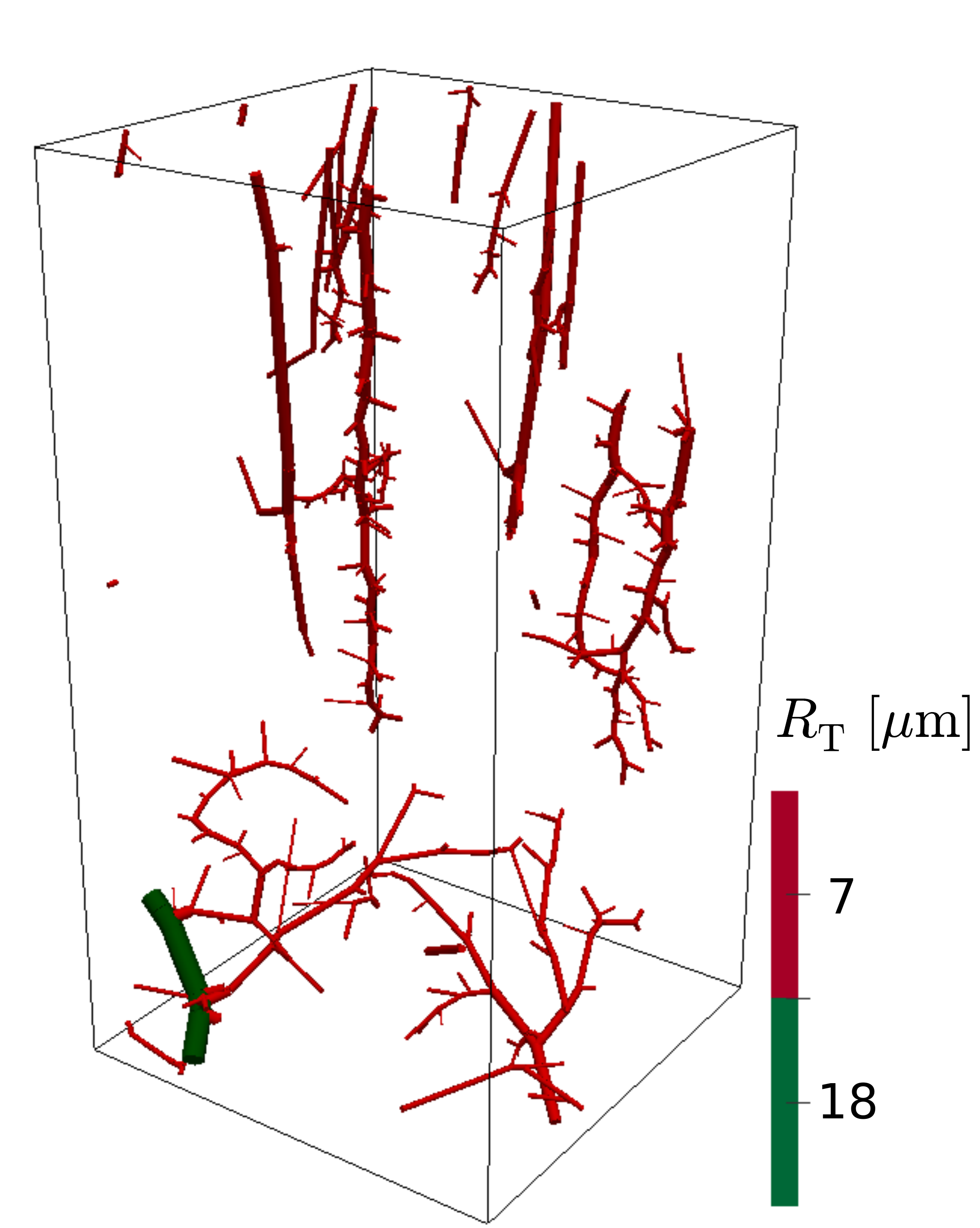}
\caption{On the left, the blood velocity for each vessel in the system $\Lambda$ is shown, where 
each segment is subject to a $1.0\;\unit{Pa}$ pressure difference. Setting the threshold to
$R_{\mathrm{T}}=7.0\;\unit{\mu m}$ yields $14918$ vessels in the set $\Lambda_{\mathrm{C}}$ with average velocity of $0.00714\;\unit{mm/s}$, 
and 337 vessels in the set $\Lambda_{\mathrm{L}}$
with an average velocity of $0.41533\;\unit{mm/s}$. On the right, the system $\Lambda_{\mathrm{L}}$ is depicted for the threshold 
$R_{\mathrm{T}}=7.0\;\unit{\mu m}$. The network
$\Lambda_{\mathrm{L}}$ for the threshold $R_{\mathrm{T}}=18.0\;\unit{\mu m}$ is marked in green. 
\label{fig:velocities}}
\end{figure}

In this subsection, we motivate the choice of the threshold $R_{\mathrm{T}}=7\;\unit{\mu m}$ that we employ to separate the larger vessels from the capillaries.
Let us consider the whole vessel system $\Lambda$ as depicted in Fig.~\ref{fig:extracted_network}, on the top. For each segment $\Lambda_k$, 
we calculate the blood velocity $v_k$, where we set a constant pressure gradient of $\delta p$ at the vertices:
$$
v_k=\frac{R_k^2}{8.0\mu_{\mathrm{bl},k}}\cdot\frac{\delta p}{|\Lambda_k|},
$$
where $\mu_{\mathrm{bl},k}$ is the viscosity of the blood according to~\eqref{eq:viscosity}. The distribution of the velocities is 
reported in Fig.~\ref{fig:velocities}, on the left,
for the case: $\delta p=1.0\;\unit{Pa}$. Choosing the threshold $R_{\mathrm{T}}$ to $7.0\;\unit{\mu m}$, we obtain that the average velocity in the set
$\Lambda_{\mathrm{C}}$ is approximately 60 times smaller than the average velocity in the set $\Lambda_{\mathrm{L}}$. Furthermore, 
the set $\Lambda_{\mathrm{C}}$ consists of
$14918$ vessels, while only $337$ vessels are contained in the set $\Lambda_{\mathrm{L}}$. Despite the low number of larger vessel, 
the chosen threshold still allows us to capture the main geometry of the vessel system, as depicted in Fig.~\ref{fig:velocities}, on the right. 
In fact, choosing a larger threshold such as $R_{\mathrm{T}}=18.0\;\unit{\mu m}$ would reduce drastically the number of larger vessels, yielding a network 
that provides only restricted informations about the geometry of the original system. On the other hand, choosing a smaller threshold would 
incorporate too many vessels in the size of capillaries. Considering the morphological values listed in \cite[Table 1]{el2018investigating}, we can observe 
that the choice of our threshold is close to the lower bound of the diameter range for the arterioles. Motivated by these considerations, we fixed the 
threshold to $7.0\;\unit{\mu m}$.

\subsection{Homogenization of the capillary bed}
\label{sec:K_up}
\begin{figure}[h!]
\centering
\includegraphics[width=0.485\textwidth]{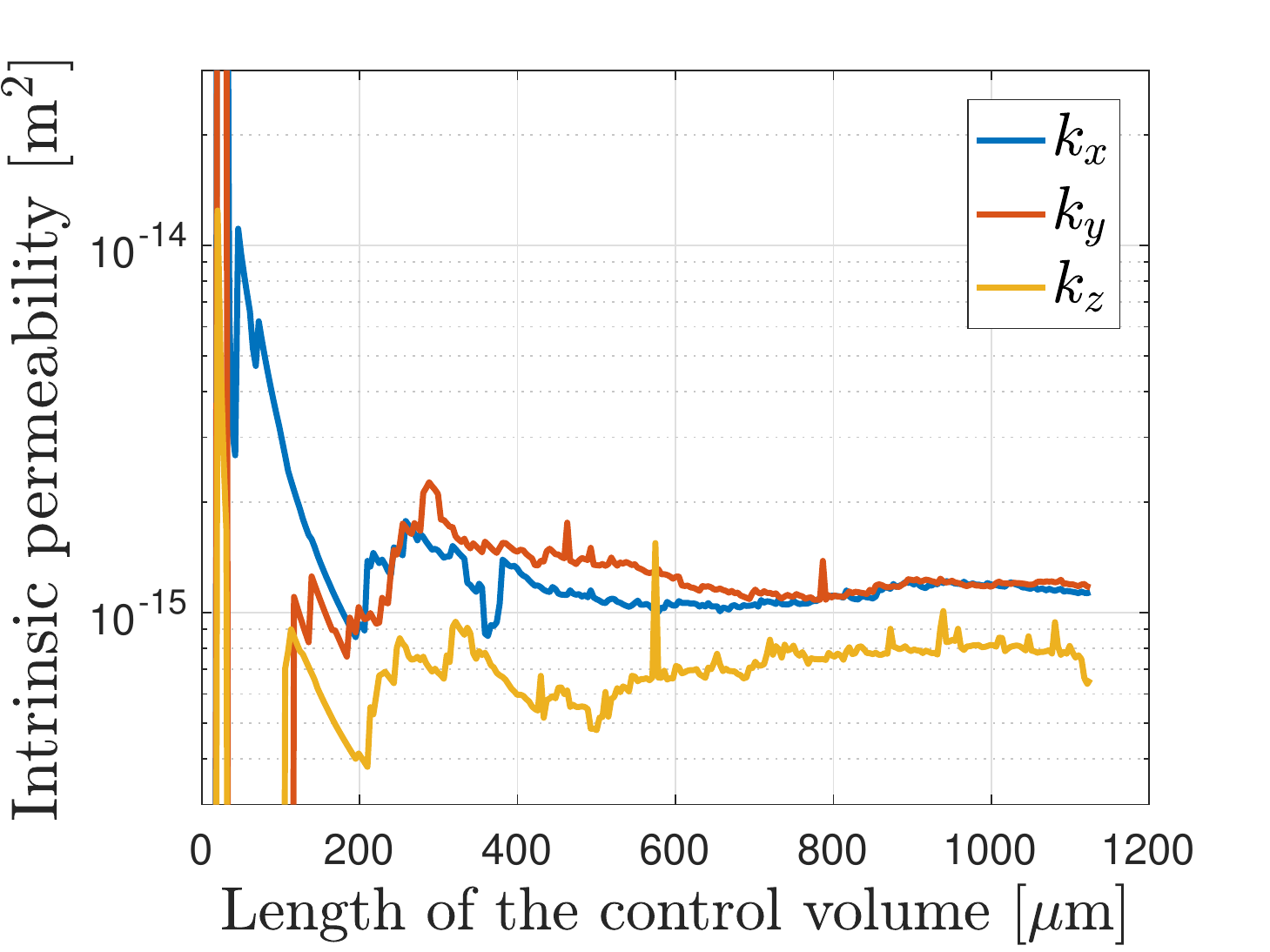}
\includegraphics[width=0.485\textwidth]{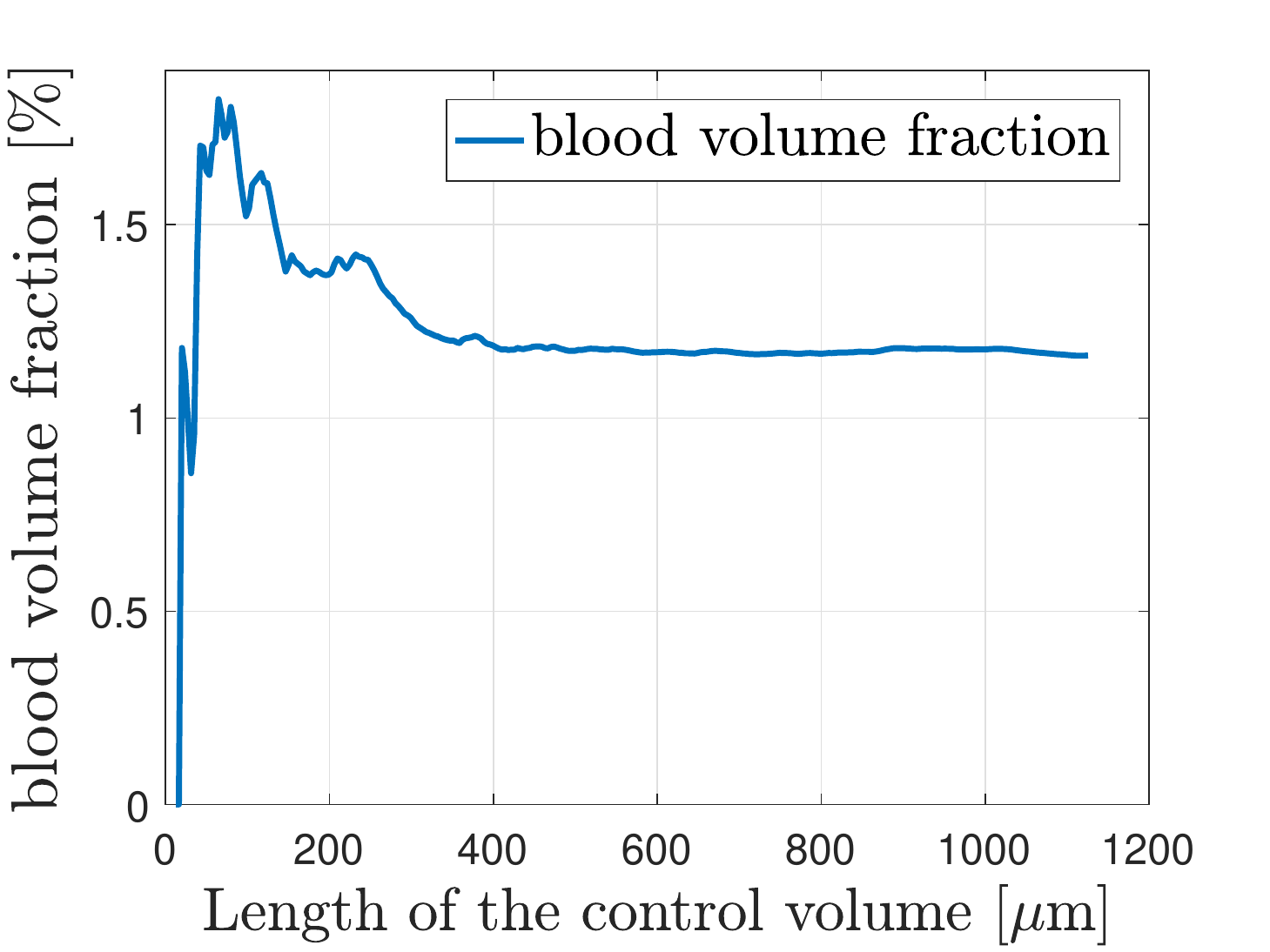}
\caption{On the left, the upscaled permeability is computed in the case where a single control volume 
is positioned in the center of the domain $\Omega$. The dimensions of the control volume
are enlarged by approximately $4\times 4\times 8\;\unit{\mu m}$ in each step until the control volume fills 
the entire domain $\Omega$. The length of the control volume in the plot has to be
doubled to obtain the actual size of the control volume in the $z$-direction. After an oscillating transition zone, 
each permeability stabilizes around a fixed value. On the right, the blood
volume fraction of the corresponding control volume is reported. \label{fig:rev_poro_perm}}
\end{figure}

\begin{figure}[!h]
\centering
\includegraphics[width=1.0\textwidth]{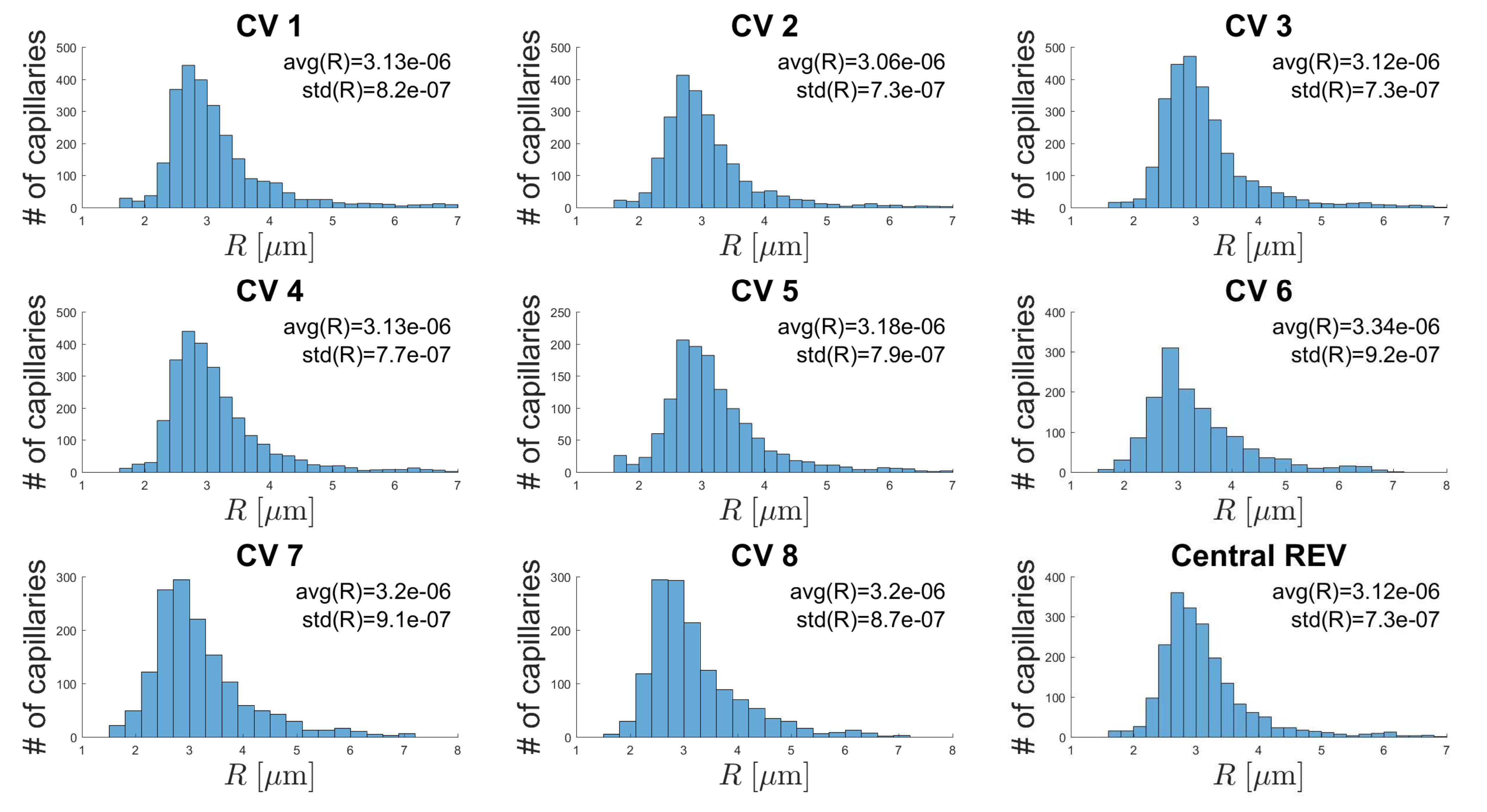}
\caption{Histrograms of the radii distributions of the capillaries contained in each control volume. The numeration is given according
to Table~\ref{tab:pressure_comp}. For each CV, the mean value avg(R) of the radii and the standard deviation std(R) are provided. 
In the last histrogram, the radii distribution of the capillaries
contained in the REV with the same center as the domain $\Omega$ and same dimensions as the other control volumes is reported.
\label{fig:radii_distributions}}
\end{figure}

\begin{figure}[!h]
\centering
\includegraphics[width=0.45\textwidth]{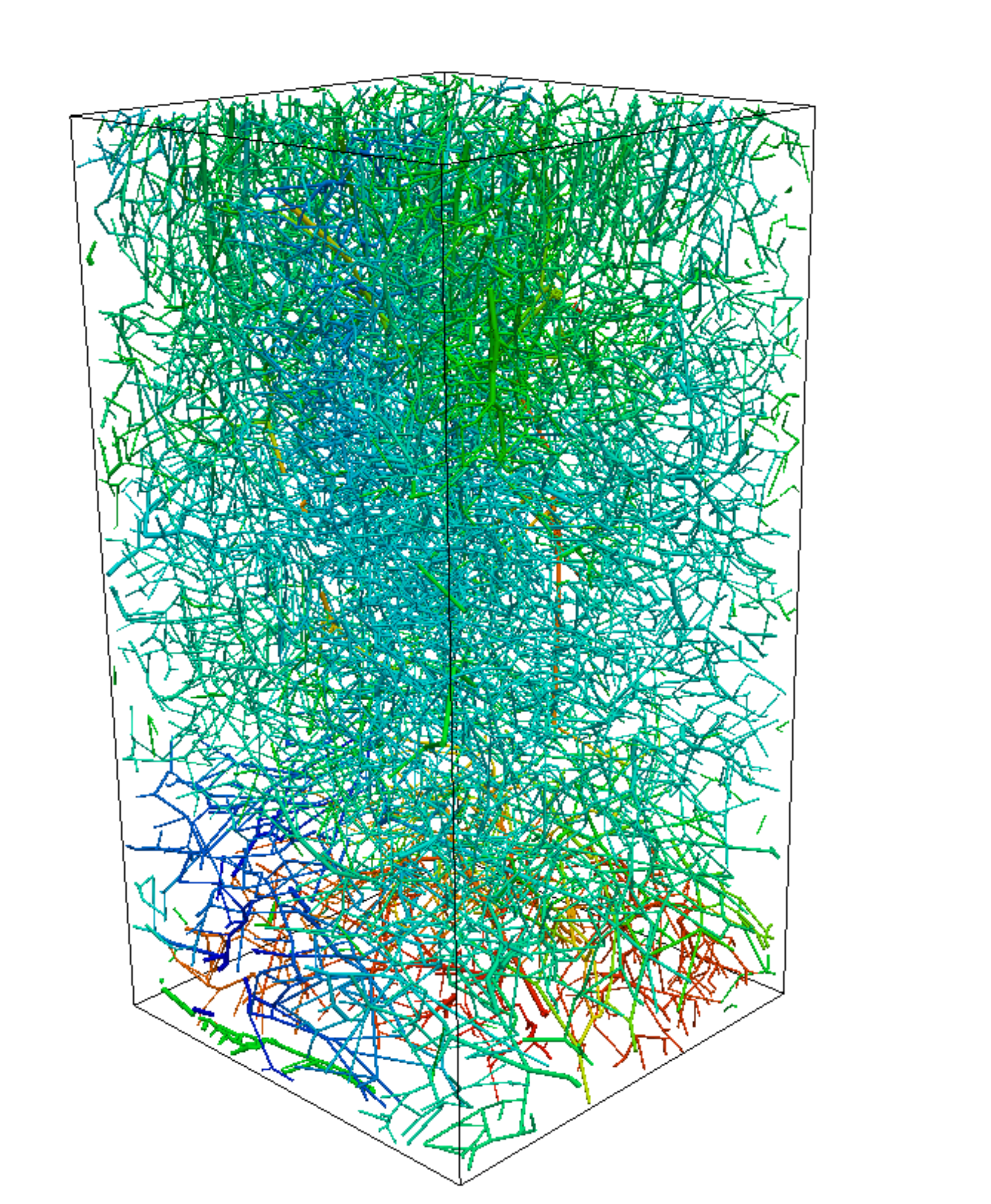}
\includegraphics[width=0.45\textwidth]{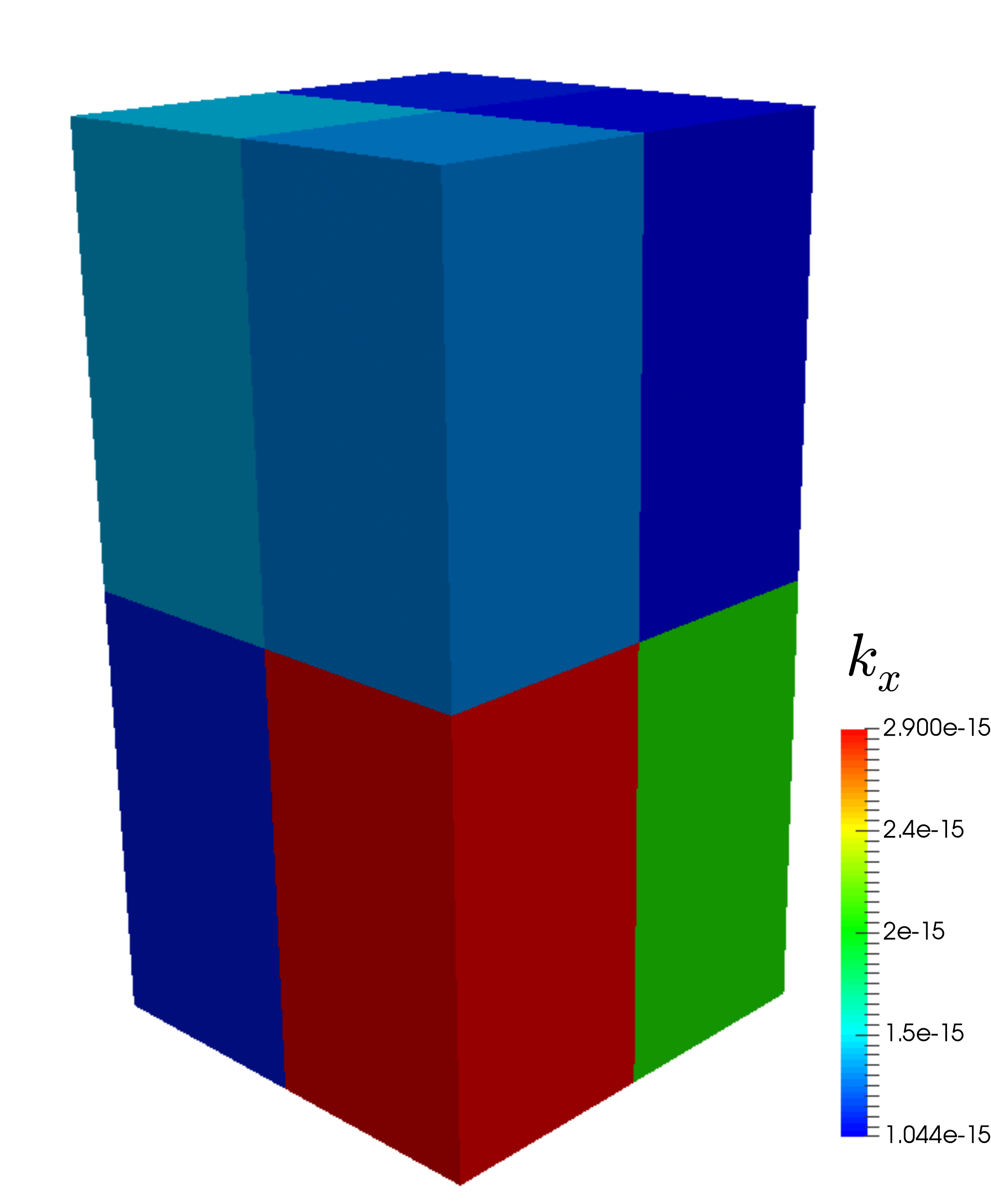}\\
\includegraphics[width=0.45\textwidth]{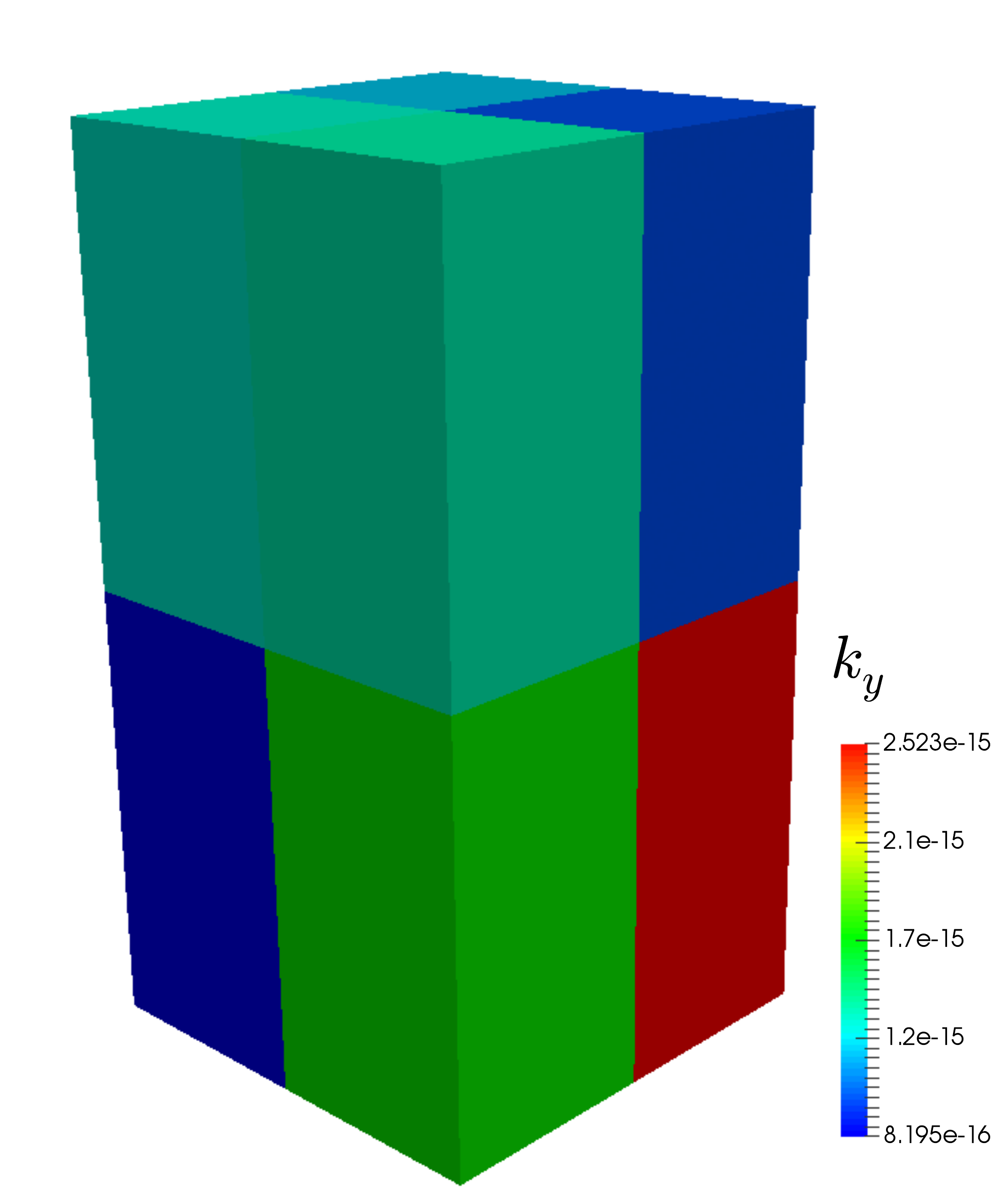}
\includegraphics[width=0.45\textwidth]{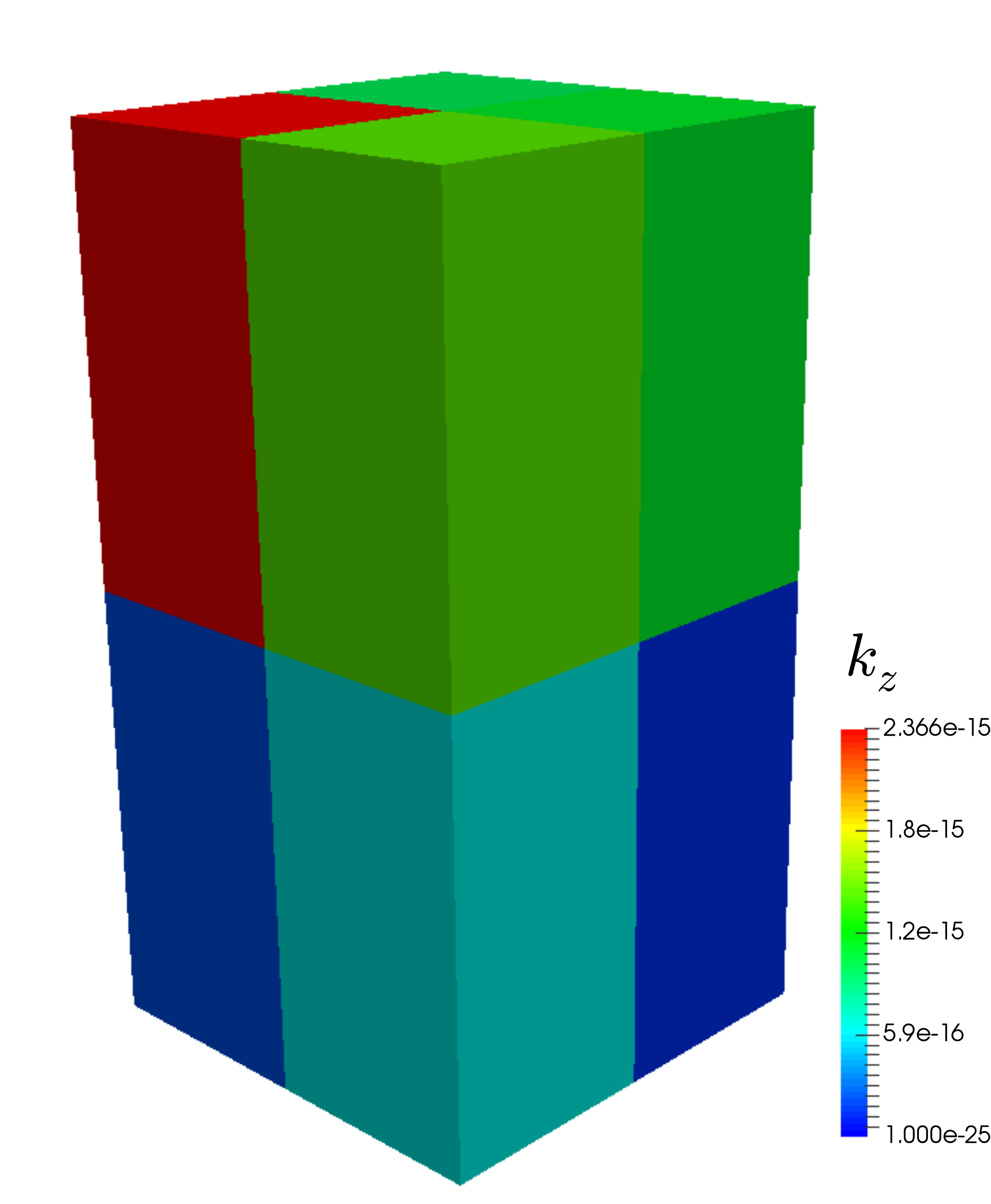}
\caption{Permeability tensors $K_{\mathrm{up}}$ in the case of $2\times 2\times 2$ REVs. The heterogeneous distribution of the capillaries
from Fig.~\ref{fig:extracted_network} is therefore incorporated, in the sense that in the lower part of the domain, the flow occurs mainly
in the $xy$-plane, while on the top in the $z$-direction.\label{fig:perm_values}}
\end{figure}

In order to determine the admissible REV sizes for the approximation of the permeability tensor, we perform the following test: A single control volume, 
initially of size
$12\times 12\times 24\;\unit{\mu m}$, is positioned in the center $(0.00056831,0.00056831,0.00113662)$ of the domain $\Omega$ and then enlarged 
in each space direction approximately by
$4.0\;\unit{\mu m}$ in the $x$- and $y$-directions and by $8.0\;\unit{\mu m}$ in the $z$-direction. For each one of these control volumes, 
the values of the intrinsic permeability are determined
using~\eqref{eq:k_up} and suitable adaptations for the $y$- and $z$-directions. In addition to that, we compute the \emph{blood volume fraction}, 
which is defined as the ratio between the
volume of the capillary network contained in the control volume under consideration and the volume of the control volume itself.
The test is performed starting from the center of $\Omega$,
because this position allows for a larger margin of growth of the control volume.
The numerical results confirm the expected oscillating behavior of the intrinsic permeability that typically occurs when the size of the 
control volume is too small (Fig.~\ref{fig:rev_poro_perm}, left).
The permeabilities appear to reach stable values, if the edges of the control volume are greater than
approximately $500.0\;\mu\unit{m}$ (in the $x$- and $y$-direction). Therefore, we can assert that the control volume with half the dimensions of the domain
($568.31 \times 568.31 \times 1081.503\;\mu \unit{m}$) and collocated in the center of $\Omega$ can be assumed to be an REV. 
A further argument to support this statement can be derived
considering the blood volume fraction of the capillary continuum. As depicted in Fig.~\ref{fig:rev_poro_perm} on the right, 
the blood volume fraction seems to stabilize around the value $1.16\%$,
if the length of the control volume is larger than approximately $400.0\;\mu\unit{m}$.
However, collocating a single REV in the center of the domain with half the sizes of $\Omega$ is not enough, since the entire 
capillary network has to be homogenized. On the other hand,
having a single REV covering the entire domain would mean that the heterogeneity of the capillary system is not considered. 
In fact, observing Fig.~\ref{fig:perm_values} (top left), we can
notice that in the upper part of the system, the capillaries are mainly aligned with the $z$-direction, while in the bottom part, 
the main directions are the $x$- and $y$-directions.
These orientations are consistent with the structure of the larger vessels, as reported in Fig.~\ref{fig:extracted_network}, on the bottom right. 
To this end, we subdivide the domain into
$2\times 2\times 2$ control volumes, each having half the sizes of the domain, as the central REV from the previous test. 
The centers of these control volumes and their corresponding numeration
are reported in Table~\ref{tab:pressure_comp}. To assert that these $8$ control volumes are REVs as well, we observe in 
Fig.~\ref{fig:radii_distributions} that the radii distribution of the
capillaries contained in each control volume is similar to that of the central REV from the previous test. Furthermore, mean radii 
and standard deviations are similar as well. Supported by
these observations, we can assume that the $8$ control volumes considered are REVs and that the permeabilities reported in 
Fig.~\ref{fig:perm_values} are representative.

\subsection{Comparison of the mass fluxes}
\label{sec:mf_comp}
A comparison between the two numerical models is provided in terms of the mass fluxes across different boundaries and interfaces.
For the one-dimensional systems in both the fully discrete network and the
reduced network, the mass flux $MF$ through a boundary node  $\mathbf{x}_k=\Lambda_k(s_k)\in \partial \Lambda\cap \partial\Omega$ is computed
as in~\eqref{eq:vgm_method}.
The inflow $MF_{\mathrm{in}}$ through the boundary point $\mathbf{x}_k$ is defined as:
$$
MF_{\mathrm{in}}(\mathbf{x}_k)=
\left\{
\begin{aligned}
&MF(\mathbf{x}_k),&\qquad&\text{if } MF(\mathbf{x}_k)>0,\\
&0,&\qquad&\text{otherwise.}\\
\end{aligned}\right.
$$
In a similar way, we can define the mass flux out of the domain:
$$
MF_{\mathrm{out}}(\mathbf{x}_k)=
\left\{
\begin{aligned}
&\left|MF(\mathbf{x}_k)\right|,&\qquad&\text{if } MF(\mathbf{x}_k)<0,\\
&0,&\qquad&\text{otherwise.}\\
\end{aligned}\right.
$$
Having this notation at hand, we define the total inflow $MF_{\mathrm{LV,in}}$ and outflow $MF_{\mathrm{LV,out}}$ through the large vessels as:
$$
\begin{aligned}
MF_{\mathrm{LV,in}}=\sum_{\mathbf{x}_k\in \partial\Lambda_\mathrm{L}\cap \partial\Omega} MF_{\mathrm{in}}(\mathbf{x}_k)\qquad\text{ and }\qquad
MF_{\mathrm{LV,out}}=\sum_{\mathbf{x}_k\in \partial\Lambda_\mathrm{L}\cap \partial\Omega} MF_{\mathrm{out}}(\mathbf{x}_k).
\end{aligned}
$$
For the hybrid approach, the mass fluxes through the boundary $\partial \Omega$ of the capillary continuum have to be interpreted as single quantities
for each boundary REV. Let us assume that the $\mathrm{REV}_j$
shares at least one side with the boundary of the domain, that is $\partial\mathrm{REV}_j\cap\partial\Omega\neq \emptyset$. In this work,
we employ the following definition for the \emph{net flux} $NF$ with respect to the $\mathrm{REV}_j$:
$$
NF(\mathrm{REV}_j)=\rho_{\mathrm{bl}}\int_{\partial\mathrm{REV}_j\cap\partial\Omega} \;
\frac{K_{\mathrm{up}}^{(j)}}{\mu_{\mathrm{bl},j}^{\mathrm{up}}}\cdot \nabla p^{\mathrm{cap}}_{(j)}\cdot \mathbf{n}\;dS,
$$
where $\mathbf{n}$ denotes the outward unit normal vector to the boundary. Numerically, the gradient $\nabla p^{\mathrm{cap}}_{(j)}$
is calculated by the standard two-point flux approximation for a
cell-centered finite volume method. Analogously as for the one-dimensional fluxes, we define the inflow $NF_{\mathrm{in}}$ through the $\mathrm{REV}_j$ as:
\begin{equation}\label{eq:nf_in}
NF_{\mathrm{in}}(\mathrm{REV}_j)=
\left\{
\begin{aligned}
&NF(\mathrm{REV}_j),&\qquad&\text{if } NF(\mathrm{REV}_j)>0,\\
&0,&\qquad&\text{otherwise,}\\
\end{aligned}\right.
\end{equation}
and the outflow as:
\begin{equation}\label{eq:nf_out}
NF_{\mathrm{out}}(\mathrm{REV}_j)=
\left\{
\begin{aligned}
&\left| NF(\mathrm{REV}_j)\right|,&\qquad&\text{if } NF(\mathrm{REV}_j)<0,\\
&0,&\qquad&\text{otherwise.}\\
\end{aligned}\right.
\end{equation}
Therefore, the total net fluxes for the capillary continuum are given by:
\begin{equation}\label{eq:cap_flux_hy}
MF_{\mathrm{cap,in}}^{\mathrm{HY}}=\sum_{j=1}^{N_{\mathrm{REV}}} NF_{\mathrm{in}}(\mathrm{REV}_j)\quad \text{and}\quad
MF_{\mathrm{cap,out}}^{\mathrm{HY}}=\sum_{j=1}^{N_{\mathrm{REV}}} NF_{\mathrm{out}}(\mathrm{REV}_j).
\end{equation}
For a suitable comparison of the fluxes, the fluxes through the capillaries in the fully-discrete method have to be averaged in 
the same sense as for the hybrid approach.
Therefore, we can similarly define the net flux $NF_{\mathrm{cap}}$ for the $\mathrm{REV}_j$ as the sum of the fluxes through the boundary
capillaries, namely:
$$
NF_{\mathrm{cap}}(\mathrm{REV}_j)=\sum_{\mathbf{x}_k\in \partial\Lambda_\mathrm{C}\cap \partial\mathrm{REV}_j\cap\partial\Omega} MF(\mathbf{x}_k).
$$
Then, the inflow $NF_{\mathrm{cap,in}}$ and outflow $NF_{\mathrm{cap,out}}$ through the $\mathrm{REV}_j$ can be defined 
analogously to~\eqref{eq:nf_in} and~\eqref{eq:nf_out}.
The total inflow $MF_{\mathrm{cap,in}}^{\mathrm{FD}}$ and outflow $MF^{\mathrm{FD}}_{\mathrm{cap,out}}$ through the capillaries 
for the fully-discrete model can be defined similarly to~\eqref{eq:cap_flux_hy}.
The mass fluxes between capillaries and tissue are denoted by $MF_{\mathrm{cap,t}}^{\mathrm{FD}}$ for the fully-discrete model and
by $MF_{\mathrm{cap,t}}^{\mathrm{HY}}$ for the hybrid model.
Again, we compute only net fluxes for each REV following a similar procedure as for the blood fluxes described above. In case of the 
hybrid model, we compute the net flux in a $\mathrm{REV}_j$ by:
$$
NF_{\mathrm{cap,t}}^{\mathrm{HY}} \left( \mathrm{REV}_j \right) = \frac{\rho_{\mathrm{int}} L_{\mathrm{cap}} S_j}
{\left| \mathrm{REV}_j \right|}\int_{ \mathrm{REV}_j } \left( p^\mathrm{t} \left( \mathbf{x} \right) -  p^\mathrm{cap} 
\left( \mathbf{x} \right) \right) - \left( \pi_p - \pi_{\mathrm{int}} \right) d\mathbf{x}.
$$
Defining the net inflow flux $NF_{\mathrm{cap,t}}^{\mathrm{HY,in}}\left( \mathrm{REV}_j \right)$ and
outflow flux $NF_{\mathrm{cap,t}}^{\mathrm{HY,out}}\left( \mathrm{REV}_j \right)$ as in the previous cases, the total inflow is given by:
\begin{equation}\label{eq:cap_t_fluxes}
MF_{\mathrm{cap,t}}^{\mathrm{HY,in}}=\sum_{j=1}^{N_{\mathrm{REV}}} NF_{\mathrm{cap,t}}^{\mathrm{HY,in}}(\mathrm{REV}_j)
\end{equation}
and the total outflow is calculated as the sum of the net outflows.\\
For the fully-discrete model, we compute the net flow in $\mathrm{REV}_j$ as follows:
$$
NF_{\mathrm{cap,t}}^{\mathrm{FD}} \left( \mathrm{REV}_j \right) = \rho_{\mathrm{int}} \cdot 2 \pi L_{\mathrm{cap}} \sum_{k \in  \Lambda_{\mathrm{C},j}}
R_k \int_{ \Lambda_{k} } \left( \overline{p^\mathrm{t}} -  p^\mathrm{v} \right) - \left( \pi_p - \pi_{\mathrm{int}} \right) dS.
$$
As in this case of the hybrid model, the total inflow is given as in~\eqref{eq:cap_t_fluxes}.
Having these definitions at hand, we first compute the mass fluxes for the FD-model. Using these values, we determine the values of 
the parameter $\alpha$ in \eqref{eq:lc_estimate} such that the
following objective functions are minimized:
$$
f_1(\alpha)=\sqrt{
\sum_{\beta}\sum_{\gamma}\frac12 \Big(MF_{\mathrm{\beta,\gamma}}^{\mathrm{HY}}-MF_{\mathrm{\beta,\gamma}}^{\mathrm{FD}}\Big)^2},
$$
for $\beta\in \{\mathrm{LV, cap}\}$ and $\;\gamma\in\{\mathrm{in, out}\}$, and
$$
f_2(\alpha)=\sqrt{\frac12 \Big(MF_{\mathrm{LV,in}}^{\mathrm{HY}}-MF_{\mathrm{LV,in}}^{\mathrm{FD}}\Big)^2}.
$$
\begin{figure}[h!]
\begin{center}
\includegraphics[width=0.59\textwidth]{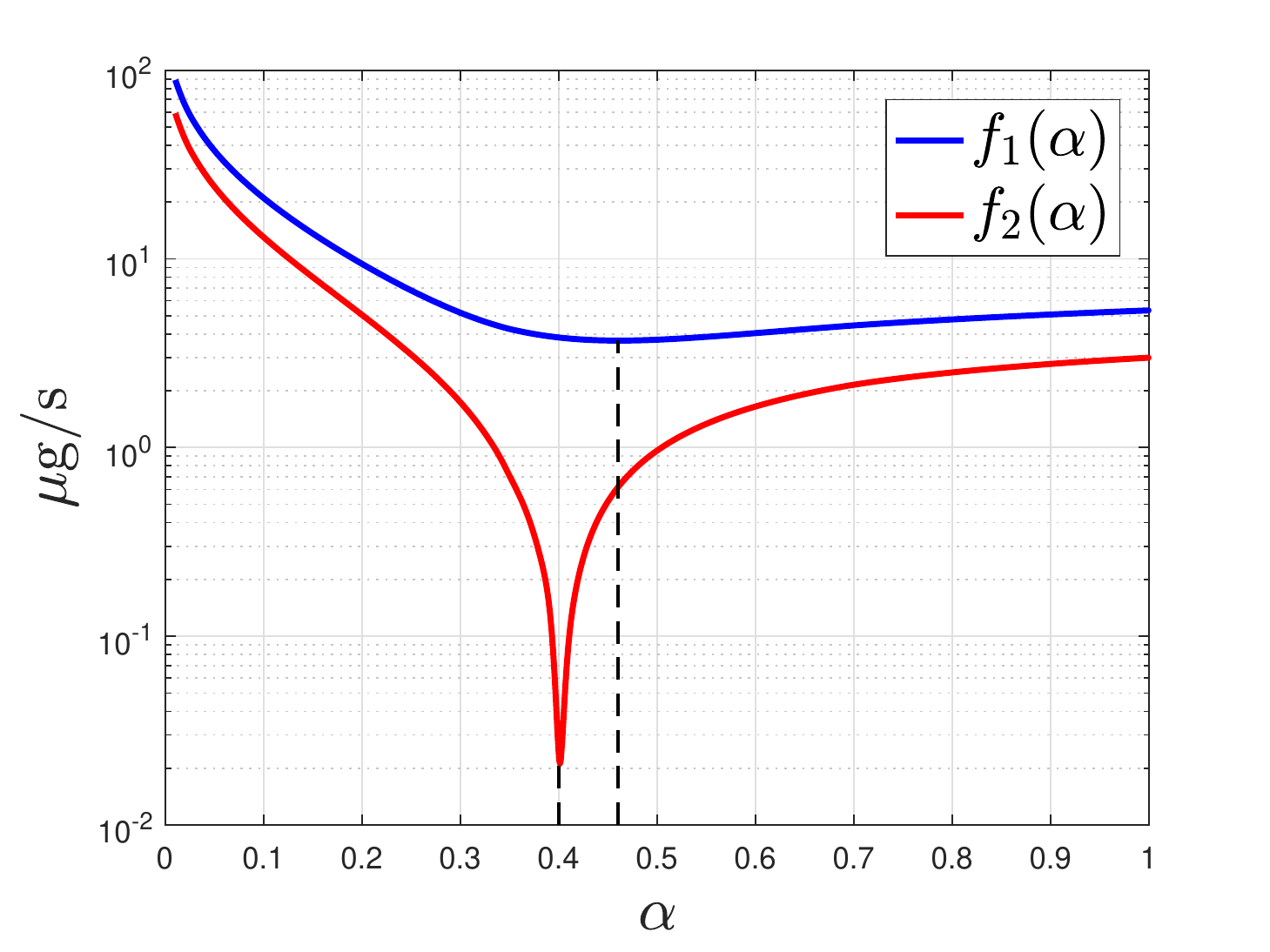}
\end{center}
\caption{\label{fig:parameter_study} Graphs of the objective functions $f_1$ and $f_2$. For $f_1$, 
the minimum is reached at $\alpha=0.46$, while, for $f_2$ it can be seen at a first glance
that the minimum is attained at $\alpha=0.4$. }
\end{figure}
The results obtained with the hybrid model strongly depend on the value of $\alpha$ in \eqref{eq:lc_estimate}, as we can 
deduce from Fig.~\ref{fig:parameter_study}, where both objective
functions $f_1$ and $f_2$ are plotted with respect to the parameter $\alpha$. For the objective function $f_2$, it is easy to 
identify the minimum at $\alpha=0.4$, while for $f_1$ the minimum
is attained at $\alpha=0.46$. The choice between these two values of $\alpha$ is made comparing the fluxes listed in 
Table \ref{table:results_hybrid_vs_fullyCoupled}, where we report the fluxes
obtained on the finest level for both the fully-discrete (FD) and the hybrid (HY) models (details to the mesh refinements are postponed to 
Section~\ref{sec:REV_pressures}). The results for the
latter one are provided for $\alpha\in\{0.2,0.4,0.46,0.9\}$. We can observe a good agreement between the hybrid and the fully-discrete models, 
in particular with respect to the inflow due to the larger vessels, if $\alpha=0.4$ is chosen. The major differences in the fluxes consist in 
the contributions of the homogenized capillaries. In fact, using the hybrid model,
the mass fluxes into the capillary continuum and out of the capillary continuum are significantly overestimated than those of the fully-discrete model. 
On the other hand, choosing $\alpha=0.46$
yields an overall better agreement with the fully-discrete model, because all four fluxes are optimized at the same time, but not a single quantity 
is as good approximated as for $\alpha=0.4$.
Moreover, due to the fact that in a small neighbourhood of the the minimum, the function $f_1$ is relatively flat, the difference between 
$\left|f_1(0.4)-f_1(0.46)\right|$ is about
$0.141\;\unit{\mu g/s}$. For completeness, we report the fluxes for $\alpha=0.2$ and $\alpha=0.9$ as well to shown by how much these results 
differ from the solution of
the fully-discrete model. Therefore, for the rest of the paper, we proceed comparing the fully-discrete model with the hybrid model, 
where we fixed $\alpha=0.4$ in~\eqref{eq:lc_estimate}.

\begin{figure}[!h]
\begin{center}
\includegraphics[width=0.465\textwidth]{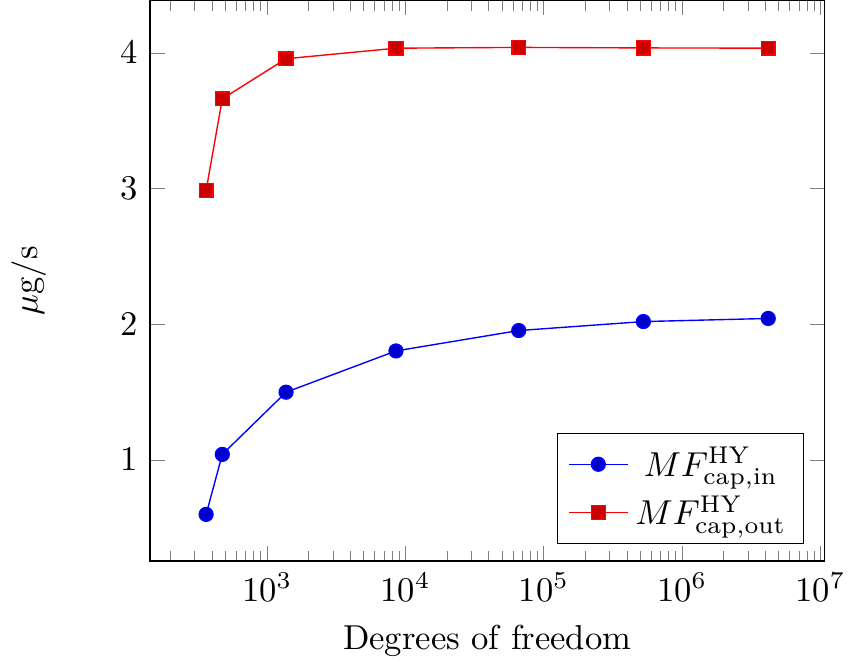}
\includegraphics[width=0.465\textwidth]{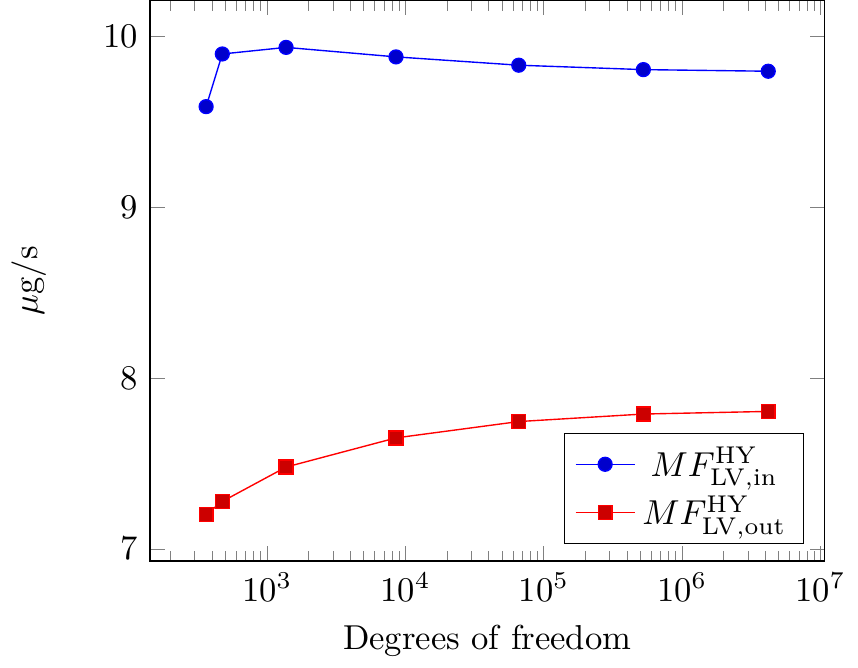}\\
\ \\
\includegraphics[width=0.465\textwidth]{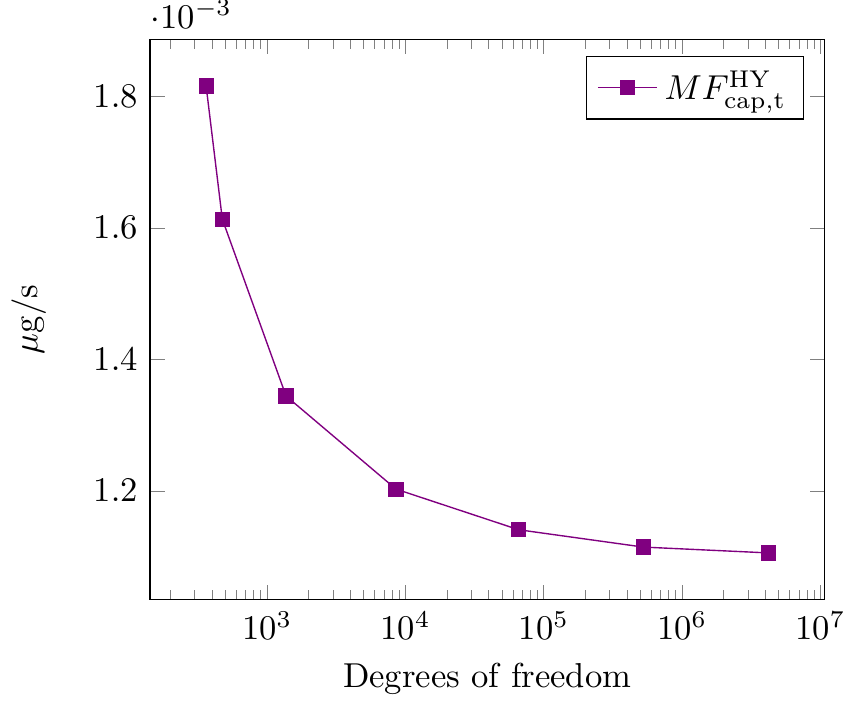}
\end{center}
\caption{Behavior of the mass fluxes for the hybrid approach with respect to the number of degrees of freedom. At the top, 
the fluxes at the boundaries of the capillaries and of the large
vessels are reported. At the bottom, the total net flux from the capillary bed into the tissue for the hybrid approach
is presented. \label{fig:convergence_fluxes}}
\end{figure}

Lastly, to validate the numerical discretizazion of the hybrid model, we report in Fig.~\ref{fig:convergence_fluxes} the numerical mass fluxes, where each plot of the flux is calculated with
respect to the mesh refinement. It can be seen that for all the quantities reported, the curves plotted are approaching asymptotic values, as the mesh is refined.
This behavior demonstrates that the mass fluxes obtained
at the finest level and reported in Table~\ref{table:results_hybrid_vs_fullyCoupled} are representative for the hybrid model.

\begin{table}[!h]
\centering
\caption{\label{table:results_hybrid_vs_fullyCoupled} Mass fluxes at the boundaries and interfaces of the vascular system. All the fluxes that are presented in this table are given
in $\unit{\mu g/s}$. For the hybrid method, the fluxes reported are obtained for different $\alpha$.}
{\small
\begin{center}
\begin{tabular}{|c||c|c||c|c||c|}
\hline
& \multicolumn{2}{|c||}{Large vessels} & \multicolumn{2}{|c||}{Capillary bed} & Tissue \\
\hline
\hline
Method & $MF_{\mathrm{LV,in}}$ & $MF_{\mathrm{LV,out}}$ & $MF_{\mathrm{cap,in}}$ & $MF_{\mathrm{cap,out}}$ & $MF_{\mathrm{cap,t}}$\\
\hline
\hline
FD     & $9.80161$ & $10.4964$ & $1.30573$  & $0.61093$ & $2.54991 \cdot 10^{-3}$\\
\hline
$\mathrm{HY}_{\alpha=0.4}$ & $9.79829$ & $7.80573$ & $2.04311$  & $4.03567$ & $1.10565 \cdot 10^{-3}$\\
\hline
$\mathrm{HY}_{\alpha=0.46}$ & $8.89951$ & $7.14353$ & $1.87204$  & $3.62801$ & $1.04549 \cdot 10^{-3}$\\
\hline
$\mathrm{HY}_{\alpha=0.2}$ & $16.9280$ & $13.2730$ & $2.86975$ & $6.52477$  & $1.48711 \cdot 10^{-3}$\\
\hline
$\mathrm{HY}_{\alpha=0.9}$ & $5.87650$ & $4.95714$ & $1.17063$ & $2.08998$  & $0.78677 \cdot 10^{-3}$\\
\hline
\end{tabular}
\end{center}
}
\end{table}

\subsection{Comparison of the REV pressures}\label{sec:REV_pressures}

After comparing the mass fluxes obtained by the two modeling approaches, we proceed with the comparison of the REV-pressures within the tissue and the capillary bed. For the hybrid model,
the averaged pressure $p_{(j),\mathrm{HY}}^ {\mathrm{cap}}$ in the capillaries for the $\mathrm{REV}_j$ is given by the definition~\eqref{eq:REV_pressure_cap}, while the average pressure
in the tissue is defined as:

\begin{equation}
\label{eq:REV_pressure_tissue}
p^{\mathrm{t}}_{(j),\mathrm{HY}} = \frac{1}{\left| \mathrm{REV}_j \right| } \int_{\mathrm{REV}_j} p^{\mathrm{t}}\left(x \right)\;dx.
\end{equation}

In case of the fully-discrete model, we use again \eqref{eq:REV_pressure_tissue} to determine the REV-pressure within the
tissue and label this value by $p_{(j),\mathrm{FD}}^{\mathrm{t}}$. The REV-pressure for the capillaries with respect to
an $\mathrm{REV}_j$ is approximated by:

$$
p_{(j),\mathrm{FD}}^{\mathrm{cap}} = \frac{1}{\left| \Lambda_{\mathrm{C},j} \right|} \int_{\Lambda_{\mathrm{C},j}} p^{\mathrm{v}}\left(x \right)\;dx.
$$

Furthermore, for the $\mathrm{REV}_j$ we define the relative pressure error $E_{\mathrm{r}}^{\mathrm{cap}}$ 
in the capillaries and $E_{\mathrm{r}}^{\mathrm{t}}$ in the tissue as:
$$
E_{\mathrm{r}}^{\beta}(j)=\frac{\left|p^{\beta}_{(j),\mathrm{HY}}-p^{\beta}_{(j),\mathrm{FD}}\right|}{p^{\beta}_{(j),\mathrm{HY}}},
\qquad \beta\in\{\mathrm{cap, t}\}.
$$
The results obtained by means of the hybrid and fully-discrete methods are reported for each REV in Table~\ref{tab:pressure_comp}, 
together with the numeration in the mesh and the
center of each REV. The average difference between the pressures obtained with the hybrid method and the fully-discrete method 
is given by approximately $537.08\;\unit{Pa}$
(corresponding to $4.03\;\unit{mmHg}$) for the capillaries, while the average difference within the tissue is given by approximately 
$316.27\;\unit{Pa}$ (corresponding to $2.37\;\unit{mmHg}$).
These values yield an average relative error of the pressures between the hybrid model and the fully-discrete of approximately 
$\overline{E^{\mathrm{cap}}_{\mathrm{r}}}=11.37\%$ in the capillaries and of $\overline{E^{\mathrm{t}}_{\mathrm{r}}}=13.97\%$ in the tissue.

\begin{table}[!h]
\centering
\caption{\label{tab:pressure_comp} Averaged REV-pressures in the capillaries and in the tissue.}
{\small
\begin{center}
\begin{tabular}{|c|c||c|c|c||c|c|c|}
\hline
\multicolumn{2}{|c||}{REV} & \multicolumn{2}{|c|}{$p^{\mathrm{cap}}_{(j)}\;[\unit{Pa}]$} & & \multicolumn{2}{|c|}{$p^{\mathrm{t}}_{(j)}\;[\unit{Pa}]$} & \\
\hline
\hline
$j$ & Center $[\unit{mm}]$  & HY & FD & $E^{\mathrm{cap}}_{\mathrm{r}}$ & HY & FD & $E^{\mathrm{t}}_{\mathrm{r}}$ \\
\hline
\hline
1 & $(0.284,0.284,0.542)$ & 5107.77 & 4535.61 & 11.20\% & 2400.12 & 1972.80 & 17.80\%\\
2 & $(0.852,0.284,0.542)$ & 5217.07 & 4704.37 & 9.83\%  & 2465.14 & 2027.71 & 17.74\%\\
3 & $(0.284,0.852,0.542)$ & 5002.77 & 4658.86 & 6.87\%  & 2325.75 & 1961.99 & 15.64\%\\
4 & $(0.852,0.852,0.542)$ & 5041.25 & 4447.64 & 11.77\% & 2392.00 & 1939.18 & 18.93\%\\
5 & $(0.284,0.284,1.623)$ & 4843.99 & 5465.34 & 12.83\% & 2156.06 & 2469.00 & 14.51\%\\
6 & $(0.852,0.284,1.623)$ & 5041.86 & 5425.96 & 7.62\%  & 2295.94 & 2440.22 & 6.28\% \\
7 & $(0.284,0.852,1.623)$ & 3789.44 & 4637.19 & 22.37\% & 1847.50 & 2192.57 & 18.68\%\\
8 & $(0.852,0.852,1.623)$ & 4960.14 & 4539.09 & 8.49\%  & 2178.97 & 2132.44 & 2.14\% \\
\hline
\end{tabular}
\end{center}
}
\end{table}

Finally, in Table~\ref{tab:pressure_comp_ref} we report the average relative errors $\overline{E^{\mathrm{cap}}_{\mathrm{r}}}$ 
and $\overline{E^{\mathrm{t}}_{\mathrm{r}}}$ of the hybrid model with mesh refinement. These errors are calculated with respect 
to the REV-pressures obtained by the fully-discrete model on the finest mesh, i.e., the values reported in Table~\ref{tab:pressure_comp} 
in the corresponding columns.
Table~\ref{tab:pressure_comp_ref} suggests that the errors $\overline{E^{\beta}_{\mathrm{r}}}$ converge to a fixed value. 
The remaining error can be considered as the modeling error arising from the homogenization.

If we calculate the solution of the hybrid model on the mesh with $16\times 16 \times 16$ elements in both the capillaries
and tissue for a total of 8538 degrees of freedom, we obtain that the average relative errors differ by less than $1\%$ from the 
average relative errors on the finest mesh. In this situation, we can assert that the modeling error dominates the discretization 
error and thus the obtained pressures can be considered as representative for the hybrid model. On the other hand, the relative error 
for the pressure in the fully-discrete model is already less than $1\%$ on the coarsest mesh, where the elements coincide with the REVs. 
However, in this case the linear system has still 12995 degrees of freedom in the network and 8 in the tissue. Therefore, compared 
with the fully-discrete model, a smaller linear system can be solved to obtain representative values for the fluxes and pressures 
in the hybrid model. This reduction in the number of degrees of freedom is expected to become more sensible, if a larger system is considered.

\begin{table}[!h]
\centering
\caption{\label{tab:pressure_comp_ref} Averaged relative errors of the REV-pressures in the capillaries and in the tissue with respect to the degrees of freedom (dofs).}
{\small
\begin{center}
\begin{tabular}{|c|c|c|}
\hline
\\[-1em]
dofs & $\overline{E^{\mathrm{cap}}_{\mathrm{r}}}$ & $\overline{E^{\mathrm{t}}_{\mathrm{r}}}$ \\
\hline
\hline
362     & 16.20\% & 22.30\% \\
\hline
474     & 14.15\% & 17.70\% \\
\hline
1370    & 12.71\% & 15.51\% \\
\hline
8538    & 11.95\% & 14.56\% \\
\hline
65882   & 11.58\% & 14.17\% \\
\hline
524634  & 11.42\% & 14.02\% \\
\hline
4194650 & 11.37\% & 13.97\% \\
\hline
\end{tabular}
\end{center}
}
\end{table}

\subsection{Sensitivity analysis of $\alpha$ for different boundary conditions}
\label{sec:boundary}
In this section, we study the influence of the boundary conditions on the parameter $\alpha$ defined in~\eqref{eq:lc_estimate}. 
Considering the results reported in Table~\ref{table:results_hybrid_vs_fullyCoupled}, one can conclude that the arterioles and venules 
determine significantly the pressure and velocity fields within the microvascular system. Therefore we vary the pressures at the 
boundary of these vessels as follows. Let us denote by $\overline{p}_{\mathrm{a}}$ and $\overline{p}_{\mathrm{v}}$ the average boundary 
pressure of the arterioles and of the venules, respectively. For the experiment setting considered here, the difference 
$$
\delta=|\overline{p}_{\mathrm{v}} -\overline{p}_{\mathrm{a}}|
$$ 
amounts to approximately $5000\;\unit{Pa}$ (corresponding to $37.5\;\unit{mmHg}$). For all $$
i\in\{-10\%,...,-1\%,0\%,1\%,...,10\%\},
$$
we add the pressure $\frac12 \delta\cdot i$ to each boundary node corresponding to an arteriole, while at 
the venous boundary nodes, we subtract the same quantity. This yields a variation in the average pressure difference by the 
fraction $i$ of $\delta$. For each new network, the optimization process described in Section~\ref{sec:mf_comp} is conducted 
and the resulting graphs of the objective functions $f_1$ and $f_2$ have the same shape as in Fig.~\ref{fig:parameter_study}. 
Following the same strategy described above, only the optimal $\alpha$ for the objective function $f_2$ is considered and reported 
in Table~\ref{tab:alpha_diff_bc}. 

If the average pressure difference is reduced, the optimal $\alpha$ is subject to relatively small variations. 
On the other hand, if $\delta$ becomes larger, the deviations from $\alpha=0.4$ become larger. However, as we can observe in 
Fig.~\ref{fig:parameter_study}, a small deviation from the optimal $\alpha$ yields a sensible difference in the flux $MF_{\mathrm{LV,in}}$. 
Therefore, a single $\alpha$ cannot be determined in advance and used for other samples with different data, but the calibration of $\alpha$ 
has to be repeated for every new experimental setting.

\begin{table}[!t]
\centering
\caption{\label{tab:alpha_diff_bc} Optimal $\alpha$ values for the cases, in which the average pressure difference $\delta$ is varied by the percentage $i$. }
{\small
\begin{center}
\begin{tabular}{|c|c|}
\hline
\\[-1em]
$i\;[\,\%\,]$ & $\alpha$\\
\hline
\hline
$-10,-9-8-7,-6,-4,-3,-2$ & 0.41\\
\hline
$-5$ & 0.42\\
\hline
$-1,0,+1$ & 0.4\\
\hline
$+2$ & 0.38\\
\hline
$+3,+4$ & 0.39\\
\hline
$+5,+6,+7$ & 0.37\\
\hline
$+8,+9,+10$ & 0.36\\
\hline
\end{tabular}
\end{center}
}
\end{table}

\section{Concluding remarks}\label{sec:Conclusion}
In this work, we have presented a hybrid model for simulating blood flow through networks at the level of microcirculation. 
The presented model is based on 1D flow models for the larger vessels and
on homogenization techniques for the capillaries and the tissue.
Thereby, the capillaries and the tissue are modeled as 3D coupled porous media resulting in a double 3D continua approach for simulating 
flows within both systems. In order to validate the simulation
results obtained by our hybrid model, we have generated a reference solution by means of a fully-discrete 3D-1D coupled model.
Here, the complete vascular network is resolved by 1D flow models and
only the flow within the tissue is considered as a porous medium flow.

For the comparison criteria between the two models, we have chosen mass fluxes at the boundaries of the microvascular system 
and averaged pressures for each REV. If the parameter $\alpha$ in~\eqref{eq:lc_estimate}
is chosen appropriately, our simulation results showed that the fluxes at the inlet and outlets of the larger vessels obtained with our hybrid 
model coincide in a satisfactory manner
with the corresponding fluxes obtained solving the fully-discrete model. On the other hand, the fluxes through 
the capillaries are overestimated by the hybrid model. Furthermore the net mass flux between the tissue and the capillaries is 
approximately $2.3$ times higher for the fully-discrete model than for the hybrid model. 
Regarding the averaged pressures for each REV, the simulations showed that the pressures
obtained with the hybrid model differ, in average, by approximately 4 mmHg in the capillaries and by approximately $2.37$ mmHg
in the tissue from the solution of the fully-discrete model. 
A more thorough comparison with respect to the reduction of the computational cost will be subject of future work. For this purpose,
a larger tissue sample should be examined. Additionally, we investigated the influence of different boundary conditions on the optimal 
parameter $\alpha$. Despite the fact that $\alpha$ varies only slightly, the results suggest that it may be necessary to calibrate $\alpha$ 
for every experimental setting to obtain accurate approximations of the mass fluxes provided by the fully-discrete model. Similarly, 
in~\cite{shipley2019hybrid} a parameter has to be optimized as well to model the flux between the different vessel types. These observations 
suggest that in context of hybrid models for microvasculature different unknown parameters occur, whose value is not known a priori. As a 
consequence, further investigations are required to improve the hybrid modeling approach. In particular, it would be of great interest to determine, 
if a combination of such parameters exists that can be applied to different settings and provides accurate results.

Considering other works regarding upscaling of capillary structures, such as~\cite{peyrounette2018multiscale}, we obtained comparable results. 
In~\cite{peyrounette2018multiscale}, the authors obtained permeabilities in the order of $10^{-14}\;\unit{m^2}$, i.e., approximately
5 times larger values than the ones depicted in Fig. 9.  
However, having a closer look to the data, one can observe that in~\cite{peyrounette2018multiscale} the radius of the capillaries 
is around $3\;\unit{\mu m}$, while in our experiment the radius of the capillaries ranges between $1.6\;\unit{\mu m}$ and $7\;\unit{\mu m}$, 
which explains the difference. Regarding the choice of the parameter $R_\mathrm{T}$, the threshold $7\;\unit{\mu m}$ is in good agreement with 
the morphological values listed in~\cite[Table 1]{el2018investigating}.

A limitation of our hybrid model consists in the determination of the parameter $\alpha$. In this paper, we employed the solution of the 
fully-discrete model to tune the parameter $\alpha$ in order to
optimize certain fluxes. A way to make a more independent definition of $\alpha$ may involve a precise computation of the quantities 
$K_{\mathrm{v}}^{(j)}$ and $\ell_{\mathrm{c}}^{(kj)}$
in~\eqref{eq:REV_pressure_cap}. Furthermore, a better approximation of the permeabilities of the homogenized system may be necessary, 
in particular, if larger systems are considered. We also point
out that the hybrid model we presented allows one to compute only net fluxes.

A clear advantage of the hybrid model is the fact that only meso-scale data are required to parametrize the model, whereas micro-scale 
data are necessary for fully discrete models. This holds for boundary data as well as for model parameters. 
Furthermore, we have provided tools to analyze homogenized models systematically that can be used to verify other upscaling strategies.

Future work in this field might be concentrated on coupling the new hybrid model for blood flow with transport models for therapeutic 
agents and other substances such that cancer therapies like hyperthermia can be simulated. A further interesting
issue that could be studied by means of the hybrid model is to enhance existing flow models for whole organs or part of organs such 
that the diagnosis techniques for clinical

\section*{Acknowledgment}
This work was partially supported by the Cluster of Excellence in Simulation Technology (EXC 310/2) and the DFG grant (WO/671 11-1). 
We gratefully acknowledge the authors of~\cite{reichold2009vascular}
and, in particular, Prof. Bruno Weber and Prof. Patrick Jenny, for providing the data set to the vascular network considered in this paper.

\bibliography{Literature}
\nocite{*}
\bibliographystyle{plain}

\end{document}